\documentclass[EJP]{ejpecp}
\pdfoutput=1

\usepackage[T1]{fontenc}
\usepackage[utf8]{inputenc}

\usepackage[british]{babel}

\usepackage{etoolbox}

\usepackage{enumerate}
\usepackage{mathrsfs}
\usepackage[retainorgcmds]{IEEEtrantools}
\usepackage{doi}

\usepackage[numbers]{natbib}

\SHORTTITLE{The hitting time of zero for a stable process}
\TITLE{The hitting time of zero for a stable process}
\AUTHORS{%
  A.\ Kuznetsov%
  \footnote{
  {\sc Department of Mathematics and Statistics, York University, Canada.}
  \BEMAIL{kuznetsov@mathstat.yorku.ca}}
  \and
  A.\ E.\ Kyprianou%
  \footnote{
    {\sc University of Bath, UK.}
    \EMAIL{a.kyprianou@bath.ac.uk}
  }
  \and
  J.\ C.\ Pardo%
  \footnote{
    {\sc  CIMAT, Mexico.}
    \EMAIL{jcpardo@cimat.mx}}
  \and
  A.\ R.\ Watson\footnote{
    {\sc University of Zurich, Switzerland.}
    \EMAIL{alexander.watson@math.uzh.ch}.
    This article was submitted while the fourth author was at the University of Bath.
  }
}

\KEYWORDS{L\'evy processes,
stable processes,
hitting times,
positive self-similar Markov processes,
Lamperti representation,
real self-similar Markov processes,
Lamperti--Kiu representation,
generalised Lamperti representation,
exponential functional,
conditioning to avoid zero}
\AMSSUBJ{60G52, 60G18, 60G51}
\SUBMITTED{March 3, 2013} 
\ACCEPTED{January 24, 2014} 

\ARXIVID{1212.5153v2}
\ARXIVPASSWORD{4t7x7}

\VOLUME{19}
\YEAR{2014}
\PAPERNUM{30}
\DOI{v19-2647}

\ABSTRACT{%
For any $\alpha$-stable L\'evy process with jumps on both sides,
where $\alpha\in(1,2)$, we find the Mellin transform of
the first hitting time of the origin and give an expression
for its density. This complements existing work
in the symmetric case and the spectrally one-sided case;
cf.\ \cite{YYY-laws, Cor-thesis} and \cite{Pes-law,Sim-hit}, respectively.
We appeal to the Lamperti--Kiu representation 
of \citet{CPR} for real-valued
self-similar Markov processes. Our main result follows by considering
a vector-valued functional equation for the Mellin transform of the
integrated exponential Markov additive process
in the Lamperti--Kiu representation. We conclude our presentation with some
applications.
}

\let\arwtemp\epsilon
\let\epsilon\varepsilon
\let\varepsilon\arwtemp

\newcommand{\RR}{\mathbb{R}}
\newcommand{\ZZ}{\mathbb{Z}}
\newcommand{\NN}{\mathbb{N}}
\newcommand{\CC}{\mathbb{C}}
\newcommand{\QQ}{\mathbb{Q}}
\newcommand{\FF}{\mathscr{F}}
\newcommand{\GG}{\mathscr{G}}
\newcommand{\Ind}{\mathbb{1}} 

\newcommand{\cdotsc}{,\dotsc,}

\newcommand{\stproca}[1]{\left(#1\right)_{t \ge 0}}
\newcommand{\stproc}[1]{\stproca{#1_t}}

\DeclareMathOperator{\rRe}{Re}
\DeclareMathOperator{\rIm}{Im}
\renewcommand{\Re}{\rRe}
\renewcommand{\Im}{\rIm}
\newcommand{\iu}{\mathrm{i}} 

\newenvironment{eqnarr}{\begin{IEEEeqnarray}{rCl}}{\end{IEEEeqnarray}\ignorespacesafterend}
\newenvironment{eqnarr*}{\begin{IEEEeqnarray*}{rCl}}{\end{IEEEeqnarray*}\ignorespacesafterend}
\newcommand{\eqnarrLHS}[1]{\IEEEeqnarraymulticol{3}{l}{#1} \\ \quad}

\renewcommand{\eqref}[1]{\hyperref[#1]{(\ref*{#1})}}

\newcommand{\xiLS}{\xi^{\mathrm{L}}}
\newcommand{\xiCPP}{\xi^{\mathrm{C}}}
\newcommand{\CE}{\Psi}
\newcommand{\CELS}{\CE^{\mathrm{L}}}
\newcommand{\CECPP}{\CE^{\mathrm{C}}}

\newcommand{\LDCPP}{\pi^{\mathrm{C}}}

\newcommand{\LSabs}{\xi^*}

\newcommand{\rhohat}{\hat{\rho}}

\newcommand{\FFt}{(\FF_t)_{t \ge 0}}
\newcommand{\GGt}{(\GG_t)_{t \ge 0}}

\newcommand{\stP}{\mathrm{P}}
\newcommand{\stPh}{\mathrm{\hat P}}
\newcommand{\stE}{\mathrm{E}}
\newcommand{\stEh}{\mathrm{\hat E}}

\newcommand{\stPaz}{\stP^\updownarrow}

\newcommand{\Xaz}{X^\updownarrow}
\newcommand{\stEaz}{\stE^\updownarrow}

\newcommand{\LevP}{\mathbb{P}}

\newcommand{\LevE}{\mathbb{E}}

\newcommand{\MAPP}{\mathbb{P}}
\newcommand{\MAPE}{\mathbb{E}}

\newcommand{\Raz}{R^\updownarrow}
\newcommand{\xiaz}{\xi^\updownarrow}

\newcommand{\CEaz}{\CE^\updownarrow}

\DeclareMathOperator{\diag}{diag}
\newcommand{\Had}{\circ}
\DeclareMathOperator{\sgn}{sgn}
\DeclareMathOperator{\Res}{Res}

\newcommand{\MT}{\mathcal{M}} 

\newcommand{\dd}{\mathrm{d}}

\newcommand{\rhoh}{\hat{\rho}}

\newcommand{\adparam}{\mathcal{A}}

\newcommand{\adparamS}{\adparam_{\text{st}}}
\newcommand*{\Indic}[1]{\Ind_{\{#1\}}}
\newcommand{\ee}{\mathbf{e}}

\newcommand{\downto}{\downarrow}

\newcommand*{\bigabs}[1]{\bigl\lvert #1 \bigr\rvert}

\newcommand{\for}{\qquad}

\newcommand{\define}{\emph}

\IEEEeqnarraydefcolsep{9}{4em}
\newcommand*{\biggabs}[1]{\biggl\lvert #1 \biggr\rvert}
\newcommand{\abs}[1]{\lvert #1 \rvert}

    \def\beq{\begin{eqnarr}}
    \def\eeq{\end{eqnarr}}
    \def\beqq{\begin{eqnarr*}} 
    \def\eeqq{\end{eqnarr*}} 
    
    \def\q{{\mathbb Q}}
    \let\re\Re
    \let\im\Im

    \def\r{{\mathbb R}}
    \def\c{{\mathbb C}}
    \def\d{{\textnormal d}}

\renewcommand{\ln}{\log}

\newcommand*{\pref}[1]{\hyperref[#1]{(\ref*{#1})}}
\newcommand*{\refpref}[2]{\hyperref[#2]{\ref*{#1}(\ref*{#2})}}

\newcommand{\arxivref}[2]{\href{http://arxiv.org/abs/#1}{arXiv:#1}%
  \ifstrempty{#2}{}{ [#2]}%
}

\let\tempFIRSTPAGE\FIRSTPAGE
\let\FIRSTPAGE\relax
\begin{document}
\tempFIRSTPAGE

\nocite{LS-hit}

\section{Introduction}
Let $X: = (X_t)_{t\geq 0}$ be a one-dimensional L\'evy process,
starting from zero,
with law $\stP$. The L\'evy--Khintchine formula
states that for all 
$\theta\in \RR$,
the characteristic exponent
$\Psi(\theta) : = 
-\log \stE (e^{\iu\theta X_1})$ satisfies 
\[
\Psi(\theta)
  = \iu a\theta 
  + \frac{1}{2}\sigma^2\theta^2 
  + \int_{\RR} (1 - e^{\iu\theta x} + \iu\theta x\Indic{|x|\leq 1})\Pi(\dd x)
  + q,
\]
where $a\in\mathbb{R}$, $\sigma\geq 0$ and $\Pi$ is a measure
(the \define{L\'evy measure}) concentrated on
$\RR\setminus\{0\}$
such that $\int_{\RR}(1\wedge x^2)\Pi(\dd x)<\infty$. The
parameter $q$ is the killing rate. When $q=0$, we say that
the process $X$ is \emph{unkilled}, and it remains in $\RR$
for all time; when $q > 0$, the process $X$ is sent to
a cemetery state at a random time, independent of the
path of $X$ and with an exponential distribution of rate $q$.

The process $(X,\stP)$ is
said to be a \define{(strictly) $\alpha$-stable process}
if it is an unkilled L\'evy process which
also satisfies the \define{scaling property}: under $\stP$,
for every $c > 0$,
the process $(cX_{t c^{-\alpha}})_{t \ge 0}$
has the same law as $X$.
It is known that $\alpha \in (0,2]$, and the case $\alpha = 2$
corresponds to Brownian motion, which we exclude.
In fact, we will assume $\alpha \in (1,2)$.
The L\'evy-Khintchine representation of such a process
is as follows:
$\sigma = 0$, and
$\Pi$ is absolutely continuous with density given by 
\[ 
  c_+ x^{-(\alpha+1)} \Indic{x > 0} + c_- \abs{x}^{-(\alpha+1)} \Indic{x < 0},
  \for x \in \RR,
\]
where $c_+,\, c_- \ge 0$, and $a = (c_+-c_-)/(\alpha-1)$.

The process $X$ has the
characteristic exponent
\begin{equation}\label{e:stable CE}
  \Psi(\theta) =
  c\abs{\theta}^\alpha
  (1  - \iu\beta\tan\tfrac{\pi\alpha}{2}\sgn\theta), \for  \theta\in\RR,
\end{equation}
where $\beta = (c_+- c_-)/(c_+ + c_-)$ and 
$c = - (c_++c_-)\Gamma(-\alpha)\cos (\pi\alpha/2)$. For more details,
see \citet[Theorems 14.10 and 14.15]{Sato}.

For consistency with the literature that we shall appeal to in this article,
we shall always parametrise our $\alpha$-stable process such that 
\[ c_+ = \Gamma(\alpha+1) \frac{\sin(\pi \alpha \rho)}{\pi} \quad \text{and} \quad
  c_- = \Gamma(\alpha+1) \frac{\sin(\pi \alpha \rhoh)}{\pi },
  \]
where
$\rho = \stP(X_t \ge 0)$ is the positivity parameter, 
and $\rhoh = 1-\rho$.  In that case, the constant $c$ simplifies to just $c = \cos (\pi\alpha(\rho - 1/2))$.

We take the point of view that the class of stable processes,
with this normalisation, is
parametrised by $\alpha$ and $\rho$; the reader will note that
all the quantities above can be written in terms of these parameters.
We shall restrict ourselves a little further within this class
by excluding the possibility of having only one-sided jumps.
This gives
us the following set of admissible parameters (see \cite[\S VII.1]{BertoinLP}):
\begin{eqnarr*}
  \adparamS &=&
  \bigl\{ (\alpha,\rho) : \alpha \in (1,2), \, \rho \in (1-1/\alpha, 1/\alpha) \bigr\}.
\end{eqnarr*}

\medskip \noindent
Let us write $\stP_x$ for the law of the shifted process $X+x$ under $\stP$.
We are interested in computing the distribution of
\[ T_0 = \inf\{ t \ge 0 : X_t = 0 \} ,\]
the first hitting time of zero for $X$, under $\stP_x$, for $x \ne 0$.
When $\alpha > 1$, this random variable is a.s.\ finite, while
when $\alpha \le 1$, points are
polar for the stable process, so $T_0 = \infty$ a.s.\ (see
\cite[Example 43.22]{Sato});
this explains our exclusion of such processes.
This paper consists of two parts: the first deals with the case
where $X$ is symmetric. Here, we may identify a positive, self-similar Markov process
$R$ and make use of the Lamperti transform to write down the
Mellin transform of $T_0$.
The second part concerns the general case where $X$ may be asymmetric.
Here we present instead
a method making use of the generalised Lamperti transform
and Markov additive processes.

It should be noted that the
symmetric case can be deduced from the general case, and so
in principle we need not go into details when $X$ is symmetric; however,
this case provides  familiar ground on which to rehearse the arguments
which appear in a more complicated situation in the general case.
Let us also note here that, in the symmetric case, the distribution of
$T_0$ has been characterised in \citet[Theorem 5.3]{YYY-laws}, and
the Mellin transform appears
in \citet[equation (1.36)]{Cor-thesis}; however, these authors proceed
via a different method.

For the spectrally one-sided case, which our range of parameters omits,
representations of law of $T_0$ have been given by \citet{Pes-law} and \citet{Sim-hit}.
This justifies
our exclusion of the one-sided case. 
Nonetheless, as we explain in Remark \ref{r:spec os}, our methodology can also
be used in this case.

We now give a short outline of the coming material.

\medskip\noindent
In section \ref{s:symm}, we suppose that the stable process $X$ is symmetric,
that is $\rho = 1/2$,
and we define a process $R$ by
\[ R_t = \abs{X_t} \Indic{t < T_0}, \for t \ge 0, \]
the \define{radial part of $X$}.
The process $R$ satisfies the $\alpha$-scaling property,
and indeed is a positive, self-similar Markov process, whose
Lamperti representation, say $\xi$, is a L\'evy process;
see section \ref{s:symm} for definitions.
It is then known that $T_0$ has the same distribution as
the random variable
\[ I(\alpha\xi) := \int_0^\infty \exp(\alpha \xi_t) \, \dd t, \]
the so-called
exponential functional of $\alpha\xi$.

In order to find the distribution of $T_0$, we compute the
Mellin transform, $\LevE[I(\alpha\xi)^{s-1}]$, for a suitable range of $s$.
The result is given in Proposition \ref{p:MT symm}.
This then characterises the
distribution, and the transform here can  be analytically inverted.

\medskip\noindent
In section \ref{s:asymm} we consider
the general case, where $X$ may not be symmetric, and our reasoning
is along very similar lines. The process $R$ still satisfies the scaling
property,
but, since its dynamics depend on the sign of $X$, it is no longer a Markov process,
and the method above breaks down. However,
due to the recent work of \citet{CPR}, there is still a type of Lamperti representation for $X$,
not in terms of  a L\'evy process, but in terms of a so-called Markov additive process, say $\xi$.
Again, the distribution of $T_0$ is equal to that of $I(\alpha\xi)$
(but now with the role of $\xi$
taken by the Markov additive process), and we develop techniques
to compute a vector-valued Mellin transform for the exponential function of
this Markov additive process.
Further, we invert the Mellin transform of $I(\alpha\xi)$ in order to deduce
explicit series representations for the law of $T_0$.

After the present article was submitted, the preprint of \citet{LS-hit} appeared,
in which the authors begin from the classical potential theoretic formula
\[ \LevE_x [e^{-q T_0}] = \frac{u^q(-x)}{u^q(0)} , \for q > 0, \; x \in \RR, \]
where $u^q$ is the $q$-potential of the stable process. Manipulating this
formula, they derive the Mellin transform of $T_0$. Their proof is rather shorter
than ours, but it appears to us that the Markov
additive point of view offers a good insight into the structure of
real self-similar Markov processes in general, and, for example,
will be central to the development of a theory of entrance laws of recurrent extensions
of rssMps.

\medskip\noindent
In certain scenarios  the distribution of $T_0$ is a very convenient quantity to have,
and we consider some applications in section \ref{s:app}:
for example, we give an alternative description of the stable process
conditioned to avoid zero, and we give some identities in law similar
to the result of \citet{BY-ent} for the entrance law
of a pssMp started at zero.

\section{The symmetric case}
\label{s:symm}

In this section, we give a brief derivation of the Mellin transform
of $T_0$ for a symmetric stable process. As we have said, we do this by
considering the Lamperti transform of the radial part of the process;
let us therefore recall some relevant definitions and results.

A \define{positive self-similar Markov process} (\define{pssMp}) with
\define{self-similarity index} $\alpha > 0$ is a standard Markov process
$R = (R_t)_{t\geq 0}$ with associated filtration $\FFt$ and probability laws
$(\stP_x)_{x > 0}$, on $[0,\infty)$, which has $0$ as an absorbing state and
which satisfies the \define{scaling property}, that for every $x, c > 0$,
\[
  \label{scaling prop}%
  \text{ the law of } (cR_{t c^{-\alpha}})_{t \ge 0}
  \text{ under } \stP_x \text{ is } \stP_{cx} \text{.}
\]
Here, we mean ``standard'' in the sense of \cite{BG-mppt},
which is to say, $\FFt$ is a complete, right-continuous filtration,
and $R$ has c\`adl\`ag paths and is strong Markov
and quasi-left-continuous.

In the seminal paper \cite{LampertiT}, Lamperti describes a one-to-one correspondence
between pssMps and L\'evy processes, which we now outline.
It may be worth noting that we have presented a slightly
different definition of pssMp from Lamperti; for the connection, see
\cite[\S 0]{VA-Ito}.

Let
$S(t) = \int_0^t (R_u)^{-\alpha}\, \dd u .$
This process is continuous and strictly increasing until $R$ reaches zero.
Let $(T(s))_{s \ge 0}$ be its inverse, and define
\[ \eta_s = \log R_{T(s)} \qquad s\geq 0. 
\]
Then $\eta : = (\eta_s)_{s\geq 0}$ is a L\'evy process started
at $\log x$, possibly killed at an independent
exponential time; the law of the L\'evy process and the
rate of killing do not depend
on the value of $x$. The real-valued process
$\eta$ with probability laws
$(\LevP_y)_{y \in \RR}$ is called the
\define{L\'evy process associated to $R$},
or the \define{Lamperti transform of $R$}.

An equivalent definition of $S$ and $T$, in terms of $\eta$ instead
of $R$, is given by taking
$T(s) = \int_0^s \exp(\alpha \eta_u)\, \dd u$
and $S$ as its inverse. Then,
\begin{equation*}
  \label{Lamp repr}
  R_t = \exp(\eta_{S(t)}) 
\end{equation*}
for all $t\geq 0$, and this shows that the Lamperti transform is a bijection.

\medskip\noindent
Let $X$ be the symmetric stable process, that is,
the process defined in the introduction with $\rho = 1/2$.
We note that this process is such that
$\stE(e^{\iu\theta X_1}) = \exp(-\abs{\theta}^{\alpha})$,
and remark briefly that it has L\'evy measure
$\Pi(\dd x) = k \abs{x}^{-\alpha-1}\,\dd x$,
where
\[ k = \Gamma(\alpha+1)\frac{\sin(\pi\alpha/2)}{\pi}
  = \frac{\Gamma(\alpha+1)}{\Gamma(\alpha/2)\Gamma(1-\alpha/2)} .\]

Connected to $X$ is an important example where the Lamperti transform
can be computed explicitly. In \citet{CC06}, the authors
compute Lamperti transforms of killed
and conditioned stable processes; the simplest of their examples,
given in \cite[Corollary 1]{CC06}, is as follows.
Let 
\[ \tau_0^- = \inf\{ t > 0 : X_t \le 0 \} , \]
\label{LSabs defined}%
and denote by $\LSabs$ the Lamperti transform of the
pssMp $\stproca{X_t \Indic{t < \tau_0^-}}$.
Then $\LSabs$ has L\'evy density
\begin{equation*}
  k e^x\abs{e^x-1}^{-(\alpha+1)}, \for x \in \RR,
\end{equation*}
and is killed at rate
$k/\alpha$. 

\medskip\noindent
Using this, we can analyse a pssMp which will give us the information we seek. Let
\[ R_t = \abs{X_t} \Indic{t < T_0}, \for t \ge 0, \]
be the radial part of the symmetric stable process. It is simple to see
that this is a pssMp. Let us denote its Lamperti transform by $\xi$.

In \citet{CPP-radial}, the authors study the L\'evy process $\xi$;
they find its characteristic function, and
Wiener-Hopf factorisation, when
$\alpha < 1$, as well as a decomposition into two simpler processes
when $\alpha \le 1$. We will now demonstrate that their expression for the
characteristic function is also valid when $\alpha > 1$, by showing
that their decomposition has meaning in terms of the Lamperti transform.


\begin{proposition}\label{p:xi decomp}
The L\'evy process $\xi$ is the sum of two independent L\'evy processes,
$\xiLS$ and $\xiCPP$, such that
\begin{enumerate}[(i)]
  \item
    The L\'evy process $\xiLS$ has characteristic exponent
    \[ \Psi^*(\theta) - k/\alpha, \for \theta \in \RR , \]
    where $\Psi^*$ is the characteristic exponent of the process
    $\xi^*$, which is the Lamperti transform of the stable process killed
    upon first passage below zero.
    That is, $\xiLS$ is formed by removing the independent
    killing from $\xi^*$.
    \item The process $\xiCPP$ is a compound Poisson process whose jumps occur
    at rate $k/\alpha$, whose L\'evy density is
    \begin{equation}\label{e:ldcpp}
      \LDCPP(y) = k \frac{e^y}{(1+e^y)^{\alpha+1}} , \for y \in \RR.
    \end{equation}
\end{enumerate}

\begin{proof}
Precisely the same argument as in \cite[Proposition 3.4]{KPW-cens}
gives the decomposition into $\xiLS$ and $\xiCPP$,
and the process $\xiLS$ is exactly as in that case.
The expression \eqref{e:ldcpp}, which determines the law of $\xiCPP$,
follows from \cite[Proposition 1]{CPP-radial},
once one observes that the computation in that article does not require
any restriction on $\alpha$.
\end{proof}

\end{proposition}

We now compute the characteristic exponent of $\xi$.
As we have mentioned, when $\alpha < 1$,
this has already been computed in \cite[Theorem 7]{CPP-radial},
but whereas in that paper the authors were concerned with computing
the Wiener-Hopf factorisation, and the characteristic function
was extracted as a consequence of this, here we
provide a proof directly from the above decomposition.

\begin{theorem}[Characteristic exponent]
The characteristic exponent of the L\'evy process $\xi$
is given by
\begin{equation}
  \CE(\theta) = 2^{\alpha}
  \frac{\Gamma(\alpha/2-\iu\theta/2)}{\Gamma(-\iu\theta/2)}
  \frac{\Gamma(1/2+\iu\theta/2)}{\Gamma((1-\alpha)/2+\iu\theta/2)},
  \for \theta \in \RR.
  \label{e:WHF low index}
\end{equation}
\begin{proof}
We consider separately the two L\'evy processes
in Proposition \ref{p:xi decomp}.

By \cite[Theorem 1]{KP-HG},
\[ \CELS(\theta) =
  \frac{\Gamma(\alpha-\iu\theta) \Gamma(1+\iu\theta)}
  {\Gamma(\alpha/2-\iu\theta)\Gamma(1-\alpha/2+\iu\theta)}
  - \frac{\Gamma(\alpha)}{\Gamma(\alpha/2)\Gamma(1-\alpha/2)}, \]
and via the beta integral,
\begin{eqnarr*}
  \CECPP(\theta)
  &=& k \int_{-\infty}^\infty (1-e^{\iu\theta y}) \pi^C( y)\dd y \\
  &=& \frac{\Gamma(\alpha+1)}{\Gamma(\alpha/2)\Gamma(1-\alpha/2)}
  \biggl[ \frac{1}{\alpha}
    - \frac{\Gamma(1+\iu\theta)\Gamma(\alpha-\iu\theta)}{\Gamma(\alpha+1)}
  \biggr] .
\end{eqnarr*}

Summing these, then using product-sum identities and \cite[8.334.2--3]{GR},
\begin{eqnarr*}
  \CE(\theta)
  &=& \Gamma(\alpha-\iu\theta)\Gamma(1+\iu\theta)
  \biggl[ \frac{\sin\bigl(\pi(\alpha/2-\iu\theta)\bigr)}{\pi}
    - \frac{\sin\bigl(\pi\alpha/2\bigr)}{\pi} \biggr] \\
  &=& \frac{2}{\pi} \Gamma(\alpha-\iu\theta)\Gamma(1+\iu\theta)
    \cos \frac{\pi(\alpha-\iu\theta)}{2} \sin \frac{-\iu\theta\pi}{2} \\
  &=& 2\pi
    \frac{\Gamma(\alpha-\iu\theta)}{\Gamma(1/2 + (\alpha-\iu\theta)/2)}
    \frac{\Gamma(1+\iu\theta)}{\Gamma(1/2 + (1+\iu\theta)/2)}
    \frac{1}{\Gamma(-\iu\theta/2)\Gamma((1-\alpha+\iu\theta)/2)}.
\end{eqnarr*}
Now, applying the Legendre--Gauss duplication  formula \cite[8.335.1]{GR} for the gamma function,
we obtain the expression in the theorem.
\end{proof}
\end{theorem}

We now characterise the law of the exponential functional
\[ I(\alpha \xi) = \int_0^\infty e^{\alpha \xi_t} \, \dd t \]
of the process $\alpha\xi$, which we do via the Mellin transform,
\[ \MT(s) = \LevE[ I(\alpha\xi)^{s - 1} ] , \]
for $s \in \CC$ whose real part lies in some open interval which we will specify.

To begin with, we observe that the Laplace exponent $\psi $ of
the process $-\alpha \xi$, that is, the function such that
$\LevE e^{-z \alpha \xi_1} = e^{\psi (z)}$,
is given by
\[ \psi (z)
  = - 2^\alpha\frac{ \Gamma(1/2-\alpha z /2) }{ \Gamma(1/2-\alpha(1+z)/2) }
  \frac{\Gamma(\alpha(1+z)/2)}{\Gamma(\alpha z/2)},
  \for \Re z \in (-1,1/\alpha).
\]
We will now proceed via the `verification result' \cite[Proposition 2]{KP-HG}:
essentially, this result says that
we must find a candidate for $\MT$ 
which satisfies the functional equation
\begin{equation}
  \MT(s+1) = -\frac{s}{\psi (-s)} \MT(s),
  \label{e:verif}
\end{equation}
for certain $s$, together with some additional conditions. Let us now state our result.

\begin{proposition}\label{p:MT symm}
The Mellin transform of $T_0$ satisfies
\begin{equation}\label{e:MT symm}
  \stE_1[T_0^{s-1}]
  = \LevE_0[I(\alpha\xi)^{s-1}]
  = \sin(\pi/\alpha)
  \frac{\cos\bigl(\frac{\pi\alpha}{2}(s-1)\bigr)}
  {\sin\Bigl(\pi\bigl(s-1+\frac{1}{\alpha}\bigr)\Bigr)}
  \frac{\Gamma(1+\alpha-\alpha s)}{\Gamma(2-s)},
\end{equation}
for $\Re s \in \bigl(-\tfrac{1}{\alpha}, 2-\tfrac{1}{\alpha}\bigr)$. 
\begin{proof}
Denote the right-hand side of \eqref{e:MT symm} by $f(s)$.
We begin by noting that $-\alpha\xi$ satisfies the Cram\'er condition
 $\psi(1/\alpha-1)=0$; the verification result
\cite[Proposition 2]{KP-HG} therefore allows us to prove the proposition
for $\Re s \in (0,2-1/\alpha)$ once we verify some conditions on
this domain. 

There are three conditions to check.
For the first, we require $f(s)$ is analytic and zero-free in
the strip $\Re s \in (0, 2-1/\alpha)$; this is straightforward. The second requires us to verify that
$f$ satisfies \eqref{e:verif}. To this end, we expand $\cos$ and $\sin$
in gamma functions (via reflection formulas \cite[8.334.2--3]{GR}) and apply the Legendre duplication  formula
 \cite[8.335.1]{GR} twice.

Finally, there is an asymptotic property which is needed.  More precisely, we need
to investigate the asymptotics of $f(s)$ as $\im s \to \infty$. To do this,
we will use the following formula
\beq\label{Gamma_asymptotics}
\lim\limits_{|y|\to \infty} \frac{|\Gamma(x+\iu y)|}{\vert y\vert^{x-\frac{1}{2}}e^{-\frac{\pi}{2}|y|}}=\sqrt{2\pi}.
\eeq
This can be derived (see \cite[Corollary 1.4.4]{AAR})
from Stirling's asymptotic formula:
\beq\label{Stirlings_formula}
\ln(\Gamma(z))=\left(z-\frac{1}{2}\right)\ln(z)-z+\frac{1}{2}\ln(2\pi)+O(z^{-1}),  
\eeq
as $ z\to \infty$ and $\abs{\arg(z)}<\pi-\delta$, for fixed $\delta > 0$.

Since Stirling's asymptotic formula is uniform in any sector $|\arg(z)|<\pi-\delta$,
it is easy to see that the 
convergence in \eqref{Gamma_asymptotics} is also uniform in $x$ belonging to a compact subset of $\r$.
Using formula \eqref{Gamma_asymptotics} we check that $|1/f(s)|=O(\exp(\pi \abs{\im s}))$ as $\im s \to \infty$,
uniformly in the strip $\Re s \in (0,2-1/\alpha)$.
This is the asymptotic result that we require.

The conditions of \cite[Proposition 2]{KP-HG} are therefore
satisfied, and it follows that the formula in the proposition
holds for $\Re s \in (0,2-1/\alpha)$. Since $f(s)$ is analytic in the wider strip 
$\Re s \in (-1/\alpha, 2-1/\alpha)$, we conclude the proof by
analytic extension.
\end{proof}
\end{proposition}

We note that the expression of \citet[equation (1.36)]{Cor-thesis} can be derived
from this result via the duplication formula for the gamma function; furthermore,
it is not difficult to deduce it from \cite[Theorem 5.3]{YYY-laws}.

Now, this Mellin transform completely characterises the law of $T_0$,
and we could at this point
invert the Mellin transform to find a series expansion  for the
density of $T_0$.
However, as we will shortly perform precisely this calculation in section
\ref{s:asymm} for the general case,
we shall leave the Mellin transform as it is and proceed
to consider what happens when $X$ may not be symmetric.

\section{The asymmetric case}
\label{s:asymm}

With the symmetric case as our model, we will now tackle the general case
where $X$ may be asymmetric. The ideas here are much the same, but the
possibility of asymmetry leads us to introduce more complicated
objects: our positive self-similar Markov processes become
real self-similar Markov processes; our L\'evy processes
become Markov additive processes; and our functional equation for
the Mellin transform (\ref{e:verif}) becomes vector-valued.

The section is laid out as follows. We
devote the first two subsections to a discussion of Markov
additive processes and their exponential functionals, and then discuss
real self-similar Markov processes and the generalised Lamperti
representation. Finally, in the last
subsection, we apply the theory which we have developed to the
problem of determining the law of $T_0$ for a general two-sided jumping stable process with $\alpha\in(1,2)$.

\subsection{Markov additive processes}
\label{ss:MAP}

Let $E$ be a finite state space and $\GGt$ a standard
filtration. 
A c\`adl\`ag process $(\xi,J)$ in $\RR \times E$
with law $\MAPP$ is called a
\define{Markov additive process (MAP)} with respect to $\GGt$
if $(J(t))_{t \ge 0}$ is a continuous-time Markov chain in $E$, and
the following property is satisfied,
for any $i \in E$, $s,t \ge 0$:
\begin{multline}\label{e:MAP}
  \text{given $\{J(t) = i\}$,
the pair $(\xi(t+s)-\xi(t), J(t+s))$ is independent of
$\GG_t$,} \\ \text{and has the same distribution as $(\xi(s)-\xi(0), J(s))$
given $\{J(0) = i\}$.}
\end{multline}

Aspects of the theory of Markov additive processes
are covered in a number of texts, among them
\cite{Asm-rp1} and \cite{Asm-apq2}.
We will mainly use the notation of
\cite{Iva-thesis},
which principally works under the assumption that
$\xi$ is spectrally negative; the results which
we quote are valid without this hypothesis, however.

Let us introduce some notation.
We write $\MAPP_i = \MAPP( \cdot \vert \xi(0) = 0, J(0) = i)$;
and if $\mu$ is a probability distribution on $E$, we write
$\MAPP_\mu = \MAPP( \cdot \vert \xi(0) = 0, J(0) \sim \mu)
  = \sum_{i \in E} \mu(i) \MAPP_i$. 
We adopt a similar convention
for expectations.

It is well-known that a Markov additive process $(\xi,J)$ also satisfies
\eqref{e:MAP} with $t$ replaced by a stopping time. Furthermore, it
has the structure given by the following proposition;
see \cite[\S XI.2a]{Asm-apq2} and \cite[Proposition 2.5]{Iva-thesis}.

\begin{proposition}
  The pair $(\xi,J)$ is a Markov additive process if and only if, for each $i,j\in E$, 
  there exist a sequence of iid L\'evy processes
  $(\xi_i^n)_{n \ge 0}$ and  a sequence of iid random variables
  $(U_{ij}^n)_{n\ge 0}$, independent
  of the chain $J$, such that if $T_0 = 0$
  and $(T_n)_{n \ge 1}$ are the
  jump times of $J$, the process $\xi$ has the representation
  \[ \xi(t) = \Indic{n > 0}( \xi(T_n -) + U_{J(T_n-), J(T_n)}^n) + \xi_{J(T_n)}^n(t-T_n),
    \for t \in [T_n, T_{n+1}),\, n \ge 0. \]
\end{proposition}
For each $i \in E$, it will be convenient to define,
on the same probability space, $\xi_i$ as a L\'evy process whose distribution is the common
  law of the $\xi_i^n$ processes in the above representation; and similarly, for each $i,j \in E$, define $U_{ij}$ to
  be a random variable having the common law of the $U_{ij}^n$ variables.

Let us now fix the following setup. Firstly, we confine ourselves
to irreducible Markov chains $J$.
Let the state space $E$ be the finite set $\{1 \cdotsc N\}$, for some $N \in \NN$.
Denote the transition rate matrix of the chain $J$ by
$Q = (q_{ij})_{i,j \in E}$.
For each $i \in E$, the Laplace exponent of the L\'evy process $\xi_i$
will be written $\psi_i$, in the sense that $e^{\psi_i(z)} = \LevE(e^{z \xi_i(1)})$, for all $z \in \CC$ for which the right-hand side exists. 
For each pair of $i,j \in E$,
define
the Laplace transform $G_{ij}(z) = \MAPE(e^{z U_{ij}})$
 of the jump
distribution $U_{ij}$,
where this exists; write $G(z)$ for the $N \times N$ matrix
whose $(i,j)$th element is $G_{ij}(z)$.
We will adopt the convention that $U_{ij} = 0$ if
$q_{ij} = 0$, $i \ne j$, and also set $U_{ii} = 0$ for each $i \in E$.

A multidimensional analogue of the Laplace exponent of a L\'evy process is
provided by the matrix-valued function
\begin{equation}\label{e:MAP F}
 F(z) = \diag( \psi_1(z) \cdotsc \psi_N(z))
  + Q \Had G(z),
\end{equation}
for all $z \in \CC$ where the elements on the right are defined,
where $\Had$ indicates elementwise multiplication, also called
Hadamard multiplication.
It is then known
that
\[ \MAPE_i( e^{z \xi(t)} ; J(t)=j) = \bigl(e^{F(z) t}\bigr)_{ij} , \for i,\,j \in E, \]
for all $z \in \CC$ where one side of the equality is defined.
For this reason, $F$ is called the \define{matrix exponent} of
the MAP $\xi$.

We now describe the existence of the \define{leading eigenvalue}
of the matrix $F$, which will play a key role in our analysis
of MAPs. This is sometimes also called the
\define{Perron–Frobenius eigenvalue};
see \cite[\S XI.2c]{Asm-apq2} and \cite[Proposition 2.12]{Iva-thesis}.

\begin{proposition}
Suppose that $z \in \CC$ is such that $F(z)$ is defined.
Then, the matrix $F(z)$ has a real simple eigenvalue $\kappa(z)$,
which is larger than the real part of all its other eigenvalues.
Furthermore, the
corresponding
right-eigenvector $v(z)$ may be chosen so that
$v_i(z) > 0$ for every $i \in E$,
and normalised such that
\begin{eqnarr}
  \pi v(z) &=& 1 \label{e:h norm}
\end{eqnarr}
where $\pi$ is the equilibrium distribution of
the chain $J$.
\end{proposition}

This leading eigenvalue features in the following probabilistic
result, which identifies a martingale (also known as the Wald martingale) and change of measure
analogous to the exponential martingale and Esscher transformation
of a L\'evy process; cf.\ \cite[Proposition XI.2.4, Theorem XIII.8.1]{Asm-apq2}.

\begin{proposition}
\label{p:mg and com}
Let
\[ M(t,\gamma) = e^{\gamma \xi(t) - \kappa(\gamma)t}
    \frac{v_{J(t)}(\gamma)}{v_{J(0)}(\gamma)} ,
  \for t \ge 0, \]
for some $\gamma$ such that the right-hand side is defined.
Then,
  $M(\cdot,\gamma)$ is a unit-mean martingale with respect to $\GGt$,
  under any initial distribution of $(\xi(0),J(0))$.
\end{proposition}

The following properties of $\kappa$ will also prove useful.

\begin{proposition}
Suppose that $F$ is defined in some open interval $D$ of $\RR$.
Then, the leading eigenvalue $\kappa$ of $F$
is smooth and convex on $D$.
\begin{proof}
Smoothness follows from results on the perturbation of eigenvalues; see
\cite[Proposition 2.13]{Iva-thesis} for a full proof. The convexity of
$\kappa$ is a consequence of the convexity properties of the entries
of $F$. The proof follows simply from \cite[Corollary 9]{BN-Perron};
see also \cite{Kin-convex,Mil-convex}.
\end{proof}
\end{proposition}

\subsection{The Mellin transform of the exponential functional}

In section \ref{s:symm}, we studied the exponential functional
of a certain L\'evy process associated to the radial part of the
stable process; now we are interested in obtaining some results
which will assist us in computing the law of an integrated exponential functional
associated to Markov additive processes.

For a MAP $\xi$, let
\[ I(-\xi) = \int_0^\infty \exp(-\xi(t)) \, \dd t. \]
One way to characterise the law of $I(-\xi)$ is via its
Mellin transform, which
we write as $\MT(s)$. This is the vector in $\RR^N$
whose $i$th element is given by
\[ \MT_i(s) = \MAPE_i[I(-\xi)^{s-1}] , \for i \in E. \]

We will shortly obtain a functional equation for $\MT$,
analogous to the functional equation \eqref{e:verif} which we saw
in section \ref{s:symm}. For L\'evy processes, proofs of the
result can be found in \citet[Proposition 3.1]{CPY97},
\citet[Lemma 2.1]{MZ-ef} and \citet[Lemma 2]{Riv-RE2};
our proof follows the
latter, making changes to account for the Markov additive
property.

We make the following assumption, which is analogous to the Cram\'er
condition for a L\'evy process; recall that $\kappa$ is the leading
eigenvalue of the matrix $F$, as discussed in section \ref{ss:MAP}.

\begin{assumption}[Cram\'er condition for a Markov additive process]
\label{a:Cramer}
There exists $z_0 < 0$ such that $F(s)$ exists on $(z_0,0)$,
and some $\theta \in (0,-z_0)$, called the \define{Cram\'er number},
such that $\kappa(-\theta) = 0$.  
\end{assumption}

Since the leading eigenvalue $\kappa$ is smooth and convex
where it is
defined, it follows also that $\kappa(-s) < 0$ for $s \in (0,\theta)$.
In particular, this renders the matrix $F(-s)$ negative definite, and hence invertible. 
Furthermore, it follows that $\kappa^\prime(0-) > 0$, and hence
(see \cite[Corollary XI.2.7]{Asm-apq2} and \cite[Lemma 2.14]{Iva-thesis}) that $\xi$
drifts to $+\infty$ independently of its initial state. This
implies that $I(-\xi)$ is an a.s.\ finite random variable.

\begin{proposition}\label{p:MT func eq}
Suppose that $\xi$ satisfies the Cram\'er condition
(Assumption \ref{a:Cramer})
with Cram\'er number $\theta\in(0,1)$. 
Then, $\MT(s)$ is finite and analytic 
when $\Re s \in (0,1+\theta)$, and we have the
following vector-valued functional equation:
\[ \MT(s+1) = - s (F(-s))^{-1} \MT(s) , \for s \in (0,\theta). \]

\begin{proof}
At the end of the proof, we shall require the existence of certain
moments of the random variable
\[ Q_t  = \int_0^t e^{- \xi(u)} \, \dd u , \]
and so we shall begin by establishing this.

Suppose that $s \in (0,\theta]$, and let $p > 1$.
Then, by the Cram\'er condition, it follows that $\kappa(-s/p) < 0$,
and hence for any $u \ge 0$, $e^{-u \kappa(-s/p)} \ge 1$.

Recall that the process
\[ M(u,z) = e^{z \xi(u) - \kappa(z)u}\frac{ v_{J(u)}(z)}{v_{J(0)}(z)} ,
  \for u \ge 0 \]
is a martingale (the Wald martingale)
under any initial distribution $(\xi(0),J(0))$, and set
\[ V(z) = \min_{j \in E} v_j(z) > 0, \]
so that for each $j \in E$, $v_j(z)/V(z) \ge 1$.

We now have everything in place to make the following calculation,
which uses the Doob maximal inequality in connection with the Wald martingale
in the third line, and the Cram\'er condition in the fourth.
\begin{eqnarr*}
\MAPE_i [ Q_t^s ]
  &\le& t^s \MAPE_i \left[ \sup_{u \le t}\bigl[ e^{-s \xi(u)/p}\bigr]^p \right] \\
  &\le& t^s \MAPE_i \left[  \sup_{u \le t}\bigl[ M(u,-s/p) v_i(-s/p)(V(-s/p))^{-1}\bigr]^p \right] \\
  &\le& t^s v_i(-s/p)^pV(-s/p)^{-p} \biggl(\frac{p}{p-1}\biggr)^p
    \MAPE_i \bigl[ M(t, -s/p)^p \bigr] \\
  &\le& t^s V(-s/p)^{-p} \biggl(\frac{p}{p-1}\biggr)^p
    e^{-tp \kappa(-s/p)} \max_{j \in J} v_j(-s/p)^p
    \MAPE_i\bigl[ e^{-s \xi(t)} \bigr]
  < \infty.
\end{eqnarr*}

Now, it is simple to show that for all $s > 0$, $t \ge 0$,
\[ \biggl( \int_0^\infty e^{-\xi(u)} \, \dd u \biggr)^s
  - \biggl( \int_t^\infty e^{-\xi(u)} \, \dd u \biggr)^s
  = s \int_0^t e^{-s \xi(u)}
    \biggl( \int_0^\infty e^{-(\xi(u+v) - \xi(u))} \, \dd v \biggr)^{s-1}
    \, \dd u . \]

For each $i \in E$, we take expectations and apply the Markov additive
property.
\begin{eqnarr*}
\eqnarrLHS{\MAPE_i \biggl[
  \biggl( \int_0^\infty e^{-\xi(u)} \, \dd u \biggr)^s
  - \biggl( \int_t^\infty e^{-\xi(u)} \, \dd u \biggr)^s \biggr]}
  &=& s \sum_{j \in E} \int_0^t
    \MAPE_i\Bigl[ e^{-s \xi(u)} ; J(u) = j \Bigr]
    \MAPE_j\biggl[ \int_0^\infty e^{-\xi(v)} \, \dd v \biggr]^{s-1}
    \, \dd u \\
  &=& s \int_0^t \sum_{j \in E} \Bigl(e^{F(-s)u} \Bigr)_{ij}
    \MAPE_j\bigl[I(-\xi)^{s-1}\bigr] \, \dd u .
\end{eqnarr*}
Since $0<s<\theta<1$, it follows that
$\bigabs{\abs{x}^s -\abs{y}^s} \leq \abs{x-y}^s$
for any $x, y\in\mathbb{R}$, and so we see that 
for each $i \in E$, the left-hand side of the above equation
is bounded by
$\MAPE_i(Q_t^s) < \infty$. Since $\bigl( e^{F(-s)u} \bigr)_{ii} \ne 0$,
it follows that $\MAPE_i[I(-\xi)^{s-1}] < \infty$ also.

If we now take $t \to \infty$, the left-hand side of the previous
equality is monotone increasing,
while on the right, the Cram\'er condition ensures that $F(-s)$
is negative definite, which is a sufficient condition for
convergence, giving the limit:
\[ \MT(s+1) = -s(F(-s))^{-1} \MT(s), \for s \in (0,\theta). \]
Furthermore, as we know the right-hand side is finite, this
functional equation allows us to conclude that $\MT(s) < \infty$
for all $s \in (0,1+\theta)$.
It then follows from the general properties of Mellin transforms that
$\MT(s)$ is finite and analytic for all $s \in \CC$ such that
$\Re s \in (0,1+\theta)$.

\end{proof}
\end{proposition}

\subsection{Real self-similar Markov processes}
\newcommand{\CPRE}{\mathcal{E}}
\newcommand{\mT}{\mathcal{T}} 

In section \ref{s:symm}, we studied a L\'evy process which was
associated through the Lamperti representation to a positive, self-similar
Markov process. Here we see that Markov additive processes also admit
an interpretation as Lamperti-type representations of \emph{real}
self-similar Markov processes.

The structure of real self-similar Markov processes has been investigated
by \citet{Chy-Lam} in the symmetric case, and \citet{CPR} in general.
Here, we give an interpretation of these authors' results
in terms of a two-state Markov additive process. We begin with some
relevant definitions, and state some of the results of these authors.

A \define{real self-similar Markov process} with \define{self-%
similarity index} $\alpha > 0$ is a standard (in the sense of \cite{BG-mppt})
Markov process $X = \stproc{X}$ with probability laws
$(\stP_x)_{x \in \RR\setminus\{0\}}$ which satisfies the
\define{scaling property}, that for all $x \in \RR \setminus \{0\}$
and $c > 0$,
\[ \text{the law of }(c X_{t c^{-\alpha}})_{t \ge 0}
  \text{ under } \stP_x \text{ is } \stP_{cx} . \]

In \cite{CPR} the authors confine their attention to processes
in `class \textbf{C.4}'. An rssMp $X$ is in this class if, for
all $x \ne 0$, $\stP_x( \exists t > 0: X_t X_{t -} < 0 ) = 1$;
that is, with probability one, the process $X$ changes
sign infinitely often.
As with the stable process, define
\[ T_0 = \inf\{t \ge 0: X_t = 0 \} .\]

Such a process may be identified, under a deformation of space and
time, with a Markov additive process which we call the
\define{Lamperti--Kiu representation} of $X$. The following
result is a simple corollary of \cite[Theorem 6]{CPR}.

\begin{proposition}
  Let $X$ be an rssMp in class \textbf{C.4} and fix $x \ne 0$.
  Define the symbol
  \[ [y] = \begin{cases}
             1, & y > 0, \\
             2, & y < 0.
           \end{cases}
  \]
  Then there exists a time-change $\sigma$, adapted to the filtration of $X$,
  such that, under the law $\stP_x$, the process 
  \[ (\xi(t),J(t)) = (\log\abs{X_{\sigma(t)}}, [X_t]) , \qquad t \ge 0, \]
  is a Markov additive process with state space $E = \{1,2\}$
  under the law $\MAPP_{[x]}$.
  Furthermore, the process $X$ under $\stP_x$ has the representation
  \[ X_t = x \exp\bigl( \xi( \tau(t))
    + \iu \pi (J( \tau(t)) + 1 )\bigr) , \for 0 \le t < T_0, \]
  where $\tau$ is the inverse of the time-change $\sigma$,
  and may be given by
\begin{equation}\label{e:Lamp time change}
  \tau(t) = \inf \biggl\{ s > 0 : \int_0^s \exp(\alpha \xi(u))
  \, \dd u > t \abs{x}^{-\alpha} \biggr\}, \for t < T_0,
\end{equation}
\end{proposition}


We observe from the expression \eqref{e:Lamp time change}
for the time-change $\tau$ that under $\stP_x$, for any $x \ne 0$,
the following identity holds for $T_0$, the hitting time of zero:
\[ \abs{x}^{\alpha} T_0 = \int_0^\infty e^{\alpha \xi(u)} \, \dd u. \]
Implicit in this statement is that the MAP on the right-hand side
has law $\mathbb{P}_1$ if $x > 0$, and law $\mathbb{P}_2$ if
$x < 0$.
This observation will be exploited in the coming
section, in which we put together the theory we have outlined
so far.

\subsection{The hitting time of zero}

We now return to the central problem of this paper: computing
the distribution of $T_0$ for a stable process. We already have in
hand the representation of $T_0$ for an rssMp as the exponential
functional of a MAP, as well as a functional equation for this
quantity which will assist us in the computation.

\medskip \noindent
Let $X$ be the stable process with parameters $(\alpha,\rho) \in \adparamS$,
defined in the introduction. We will restrict our attention for now
to $X$ under the measures $\stP_{\pm 1}$; the results for other initial
values can be derived via scaling.

Since $X$ is an rssMp, it has a representation in terms of a MAP
$(\xi,J)$; futhermore, under $\stP_{\pm 1}$,
\[ T_0 = \int_0^\infty e^{\alpha \xi(s)} \, \dd s
  = I(\alpha \xi);
\]
to be precise, under $\stP_1$ the process $\xi$ is under $\MAPP_1$,
while under $\stP_{-1}$ it is under $\MAPP_2$.

In \cite[\S 4.1]{CPR}, the authors calculate the characteristics
of the Lamperti--Kiu representation for $X$, that is, the processes 
$\xi_i$, and the jump distributions $U_{ij}$ and rates $q_{ij}$.
Using this information, and the representation \eqref{e:MAP F},
one sees that the MAP $(-\alpha \xi,J)$ has matrix exponent
\[
  F(z) =
  \begin{pmatrix}
    - \dfrac{\Gamma(\alpha(1+z))\Gamma(1-\alpha z)}
      {\Gamma(\alpha\rhoh+\alpha z)\Gamma(1-\alpha\rhoh- \alpha z)}
    & \dfrac{\Gamma(\alpha(1+z))\Gamma(1-\alpha z)}
      {\Gamma(\alpha\rhoh)\Gamma(1-\alpha\rhoh)}
    \\
    \dfrac{\Gamma(\alpha(1+z))\Gamma(1-\alpha z)}
      {\Gamma(\alpha\rho)\Gamma(1-\alpha\rho)}
    & - \dfrac{\Gamma(\alpha(1+z))\Gamma(1-\alpha z)}
      {\Gamma(\alpha\rho+\alpha z)\Gamma(1-\alpha\rho-\alpha z)}
  \end{pmatrix} ,
\]
for $\Re z \in (-1,1/\alpha)$.

\begin{remark}\label{r:spec os}
  It is well-known that, when $X$ does not have one-sided jumps,
  it changes sign infinitely often; that is, the rssMp $X$ is in \cite{CPR}'s
  class \textbf{C.4}. When the stable process has only one-sided jumps, which
  corresponds to the parameter values $\rho = 1-1/\alpha,\,1/\alpha$, 
  then it jumps over $0$ at most once before hitting it; the
  rssMp is therefore in class \textbf{C.1} or \textbf{C.2} according to the
  classification of \cite{CPR}.
  The Markov chain component of
  the corresponding MAP then has one absorbing state,
  and hence is no longer irreducible.
  Although it seems plain that our calculations can be carried
  out in this case, we omit it for the sake of simplicity.
  As we remarked in the introduction, it is considered in
  \cite{Pes-law,Sim-hit}.
\end{remark}

\smallskip \noindent
We now analyse $F$ in order to deduce the Mellin transform of $T_0$.
The equation $\det F(z) = 0$ is equivalent to
\[
  \sin(\pi(\alpha\rho+\alpha z))\sin(\pi(\alpha\hat\rho+\alpha z)) 
  - \sin(\pi \alpha \rho) \sin(\pi\alpha\hat\rho)
  = 0,
\]
and considering the solutions of this, it is not difficult to deduce that
$\kappa(1/\alpha-1) = 0$; that is, $-\alpha\xi$ satisfies the
Cram\'er condition (Assumption \ref{a:Cramer}) with
Cram\'er number $\theta = 1-1/\alpha$.

Define
\[
f_1(s) := \stE_1[ T_0^{s-1} ] = \MAPE_1[ I(\alpha\xi)^{s-1}],
\quad
f_2(s) := \stE_{-1} [ T_0^{s-1}] = \MAPE_2[ I(\alpha\xi)^{s-1}],
\]
which by Proposition \ref{p:MT func eq} are defined when $\Re s \in (0,2-1/\alpha)$.
This proposition also implies that 
\beq\label{matrix_equation2}
{\mathbf B}(s)\times \left[ \begin{array}{c} f_1(s+1) \\ f_2(s+1)  \end{array} \right]=
\left[ \begin{array}{c} f_1(s) \\ f_2(s)  \end{array} \right],
\for s \in (0,1-1/\alpha),
\eeq
where $\mathbf{B}(s):=-F(-s)/s$. Using the reflection formula for the gamma function we find that 
\beqq
{\mathbf B}(s)=\frac{\alpha}{\pi} \Gamma(\alpha-\alpha s) \Gamma(\alpha s)
\left[ \begin{array}{cc} \sin(\pi \alpha (\hat \rho-s)) & -\sin(\pi \alpha \hat \rho) \\
-\sin(\pi \alpha \rho)   & \sin(\pi \alpha ( \rho-s)) \end{array} \right],
\eeqq
for $\Re s \in (-1/\alpha,1), \, s \ne 0$, and
\beq\label{det_Bs}
{\textnormal{det}}({\mathbf B}(s))=-\alpha^2 \frac{\Gamma(\alpha-\alpha s) \Gamma(\alpha s)}
{\Gamma(1-\alpha+\alpha s) \Gamma(1-\alpha s)}, \for \Re s \in (-1/\alpha,1), \, s \ne 0.
\eeq
Therefore, if we define $\mathbf{A}(s) = (\mathbf{B}(s))^{-1}$, we have
\[
{\mathbf A}(s)
= -\frac{1}{\pi \alpha} \Gamma(1-\alpha+\alpha s) \Gamma(1-\alpha s)
\begin{bmatrix} \sin(\pi \alpha (\rho-s)) & \sin(\pi \alpha \hat \rho) \\
\sin(\pi \alpha \rho)   & \sin(\pi \alpha (\hat \rho-s)) \end{bmatrix}
\]
for $\Re s \in (1-2/\alpha,1-1/\alpha)$, and may rewrite \eqref{matrix_equation2} in the form
\beq\label{matrix_equation}
\left[ \begin{array}{c} f_1(s+1) \\ f_2(s+1)  \end{array} \right]={\mathbf A}(s) \times 
\left[ \begin{array}{c} f_1(s) \\ f_2(s)  \end{array} \right],
\for s \in (0,1-1/\alpha).
\eeq

The following theorem is our main result. 
\begin{theorem}\label{thm_Mellin_transform}
For $-1/\alpha<\re(s)<2-1/\alpha$ we have
\beq\label{Mellin_transform}
  \stE_1[T_0^{s-1}]
  =  \frac{\sin\left(\frac{\pi}{\alpha} \right)}{\sin(\pi \hat \rho)}
  \frac{\sin\left( \pi \hat \rho (1-\alpha+\alpha s)\right)}
    {\sin\left(\frac{\pi}{\alpha} (1-\alpha+\alpha s) \right)}
  \frac{\Gamma(1+\alpha-\alpha s)}{\Gamma(2-s)}.
\eeq
The corresponding expression for $\stE_{-1}[T_0^{s-1}]$ can be obtained from \eqref{Mellin_transform}
by changing $\hat \rho \mapsto \rho$.
\end{theorem}

Let us denote the function in the right-hand side of \eqref{Mellin_transform} by $h_1(s)$,
and by $h_2(s)$
the function obtained from $h_1(s)$ by replacing $\hat \rho \mapsto \rho$.
Before we are able to prove Theorem 
\ref{thm_Mellin_transform}, we need to establish several properties of these functions. 

\begin{lemma}\label{lemma1} 
\mbox{}
\begin{itemize}
\item[(i)] There exists $\epsilon \in (0,1-1/\alpha)$ such that the functions $h_1(s)$, $h_2(s)$ are analytic and zero-free in the vertical strip $0 \le \re(s) \le 1+\epsilon$. 
\item[(ii)] For any $-\infty< a<b < +\infty$ there exists $C>0$ such that 
\beqq
e^{-\pi |\im(s)|}<|h_i(s)|< e^{-\frac{\pi}{2}(\alpha-1) |\im(s)|}, \;\;\; i=1,2
\eeqq
 for all $s$ in the vertical strip $a\le \re(s) \le b$ satisfying $|\im(s)|>C$.
\end{itemize}
\end{lemma}
\begin{proof}
It is clear from the definition of $h_1(s)$ that it is a meromorphic function.
Its zeroes are 
contained in the set 
\beqq
\{2,3,4,\dots\}\cup \{1-1/\alpha+n/(\alpha \hat \rho) : \; n \in {\mathbb Z}, n\ne 0 \}
\eeqq 
and its poles are contained in the set
\beqq
\{1+n/\alpha : \; n\ge 1\} \cup \{n-1/\alpha : \; n \in {\mathbb Z}, n \ne 1\}. 
\eeqq
In particular, $h_1(s)$ possesses neither zeroes nor poles in the strip $0 \le \re(s) \le 1$. 
The same is clearly true for $h_2(s)$, which proves part (i). 

We now make use of Stirling's formula, as we did in section \ref{s:symm}.
Applying \eqref{Gamma_asymptotics} to $h_1(s)$ we find 
that as $s \to \infty$ (uniformly in the strip $a\le \re(s) \le b$) we have
\beqq
\ln (|h_1(s)|)=- \frac{\pi}{2} (1+\alpha - 2\alpha \hat \rho) |\im(s)|+O(\ln(|\im(s)|)).
\eeqq
Since for $\alpha>1$, the admissible parameters $\alpha$, $\rho$ must satisfy $\alpha-1<\alpha \hat \rho<1$, this shows that
\beqq
\alpha-1<1+\alpha - 2\alpha \hat \rho<3-\alpha<2,
\eeqq
and completes the proof of part (ii). 
\end{proof}

\begin{lemma}\label{lemma2}
The functions $h_1(s)$, $h_2(s)$ satisfy the system of equations  \eqref{matrix_equation}. 
\end{lemma}
\begin{proof}
Denote the elements of the matrix $\mathbf{A}(s)$ by $A_{ij}(s)$.
Multiplying the first row 
of $\mathbf{A}(s)$ by the column vector $[h_1(s), h_2(s)]^T$, and using identity $\sin(\pi \rho)=\sin(\pi \hat \rho)$, we obtain
\begin{eqnarr*}
\eqnarrLHS{A_{11}(s)h_1(s)+A_{12}(s)h_2(s)}
  &=&
  -\frac{1}{\pi \alpha} \frac{\sin\left(\frac{\pi}{\alpha} \right)}{\sin(\pi \hat \rho)}
  \frac{\Gamma(1-\alpha s)}{\sin\left(\frac{\pi}{\alpha} (1-\alpha+\alpha s) \right)}
  \left[\frac{\Gamma(1+\alpha-\alpha s)}{\Gamma(2-s)}
  \Gamma(1-\alpha+\alpha s)  \right] \\
  && {} \times 
  \big\{ \sin(\pi \alpha (\rho-s)) \sin\left( \pi \hat \rho (1-\alpha+\alpha s)\right)
  +\sin(\pi \alpha \hat \rho) \sin\left( \pi  \rho (1-\alpha+\alpha s)\right)\big\}.
\end{eqnarr*}
Applying identity $\Gamma(z+1)=z \Gamma(z)$ and reflection
formula for the gamma function, we rewrite the expression in the square brackets as follows:
\beqq
\left[\frac{\Gamma(1+\alpha-\alpha s)}{\Gamma(2-s)}
\Gamma(1-\alpha+\alpha s)  \right] = 
\frac{\alpha \Gamma(\alpha-\alpha s)}{\Gamma(1-s)}
\Gamma(1-\alpha+\alpha s)=\frac{\pi \alpha}{\sin(\pi \alpha(1-s)) \Gamma(1-s)}.
\eeqq
Applying certain trigonometric identities, we obtain
\begin{multline*}
\sin(\pi \alpha (\rho-s)) \sin\left( \pi \hat \rho (1-\alpha+\alpha s)\right)+
\sin(\pi \alpha \hat \rho) \sin\left( \pi  \rho (1-\alpha+\alpha s)\right) \\
=\sin(\pi \alpha (1-s)) \sin(\pi \hat \rho (1+\alpha s)). 
\end{multline*}
Combining the above three formulas we conclude 
\beqq
A_{11}(s)h_1(s)+A_{12}(s)h_2(s)=
-\frac{\sin\left(\frac{\pi}{\alpha} \right)}{\sin(\pi \hat \rho)}
  \frac{\sin\left( \pi \hat \rho (1+\alpha s)\right)}
    {\sin\left(\frac{\pi}{\alpha} (1-\alpha+\alpha s) \right)}
  \frac{\Gamma(1-\alpha s)}{\Gamma(1-s)}=h_1(s+1). 
\eeqq
The derivation of the identity $A_{21}(s)h_1(s)+A_{22}(s)h_2(s)=h_2(s+1)$ is identical.
We have now established that two functions $h_i(s)$ satisfy
the system of equations \eqref{matrix_equation}.
\end{proof} 

\begin{proof}[Proof of Theorem \ref{thm_Mellin_transform}]
Our goal now is to establish the uniqueness of solutions to the system
\eqref{matrix_equation} in a certain class of meromorphic functions,
which contains both $h_i(s)$ and $f_i(s)$.
This will imply  $h_i(s) \equiv f_i(s)$. 
Our argument is similar in spirit to the proof of Proposition 2 in \cite{KP-HG}.

First of all, we check that there exists $\epsilon\in (0,1/2-1/(2\alpha))$, such that
 the functions $f_1(s)$, $f_2(s)$ are analytic and bounded in the open strip 
 \beqq
 {\mathcal S}_{\epsilon}=\{s\in \c: \epsilon < \Re(s) < 1+2\epsilon\}
 \eeqq
This follows from Proposition \ref{p:MT func eq} and the estimate
\[
\abs{f_1(s)}=\abs{\stE_1[ T_0^{s-1} ]}\le  \stE_1[ \abs{T_0^{s-1}} ]=\stE_1[ T_0^{\re(s)-1} ]=f_1(\Re(s)).
\]
The same applies to $f_2$. Given results of Lemma \ref{lemma1},
we can also assume that $\epsilon$ is small enough, so that the functions $h_i(s)$ 
are analytic, zero-free and bounded in the strip ${\mathcal S}_{\epsilon}$.

Let us define $D(s):=f_1(s)h_2(s) - f_2(s) h_1(s)$ for $s\in {\mathcal S}_{\epsilon}$. 
From the above properties of $f_i(s)$ and $h_i(s)$ we conclude that $D(s)$ is analytic and bounded in ${\mathcal S}_{\epsilon}$. Our first goal is to show that $D(s)\equiv 0$. 

If both $s$ and $s+1$ belong to ${\mathcal S}_{\epsilon}$, then the function $D(s)$ satisfies the equation
\beq\label{eqn_Fs}
D(s+1)=-\frac{1}{\alpha^2} \frac{\Gamma(1-\alpha+\alpha s) \Gamma(1-\alpha s)}{\Gamma(\alpha-\alpha s) \Gamma(\alpha s)}D(s),
\eeq
as is easily established by taking determinants of the matrix equation  
\beqq
\left[ \begin{array}{cc} f_1(s+1) & h_1(s+1) \\ f_2(s+1) & h_2(s+1)  \end{array} \right]={\mathbf A}(s) \times 
\left[ \begin{array}{cc} f_1(s) & h_1(s) \\ f_2(s) & h_2(s)  \end{array} \right], 
\eeqq
and using \eqref{det_Bs} and the identity $\mathbf{A}(s) = \mathbf{B}(s)^{-1}$.

Define also
\[
  G(s):=\frac{\Gamma(s-1)\Gamma(\alpha-\alpha s)}{\Gamma(1-s)\Gamma(-\alpha + \alpha s)}
  \sin\left( \pi \left(s+\frac{1}{\alpha} \right) \right).
\]
It is simple to check that:
\begin{enumerate}[(i)]
\item\label{G1} $G$ satisfies the functional equation \eqref{eqn_Fs};
\item\label{G2} $G$ is analytic and zero-free in the strip ${\mathcal S}_{\epsilon}$;
\item\label{G3} $|G(s)|\to \infty$ as $\im(s)\to \infty$, uniformly in the strip  ${\mathcal S}_{\epsilon}$
(use \eqref{Gamma_asymptotics} and $\alpha > 1$). 
\end{enumerate}

We will now take the ratio of $D$ and $G$ in order to obtain a periodic function,
borrowing a technique from the theory of functional equations
(for a similar idea applied to the gamma function, see \cite[\S 6]{Art-gamma}.)
We thus define $H(s):=D(s)/G(s)$ for $s\in {\mathcal S}_{\epsilon}$.
The property \pref{G2} guarantees that $H$ is analytic in the strip $ {\mathcal S}_{\epsilon}$, 
while property \pref{G1} and \eqref{eqn_Fs} show that\linebreak 
$H(s+1)=H(s)$ 
if both $s$ and $s+1$ belong to $ {\mathcal S}_{\epsilon}$.
Therefore, we can extend $H(s)$ to an entire function satisfying
$H(s+1)=H(s)$ for all $s\in \c$.
Using the periodicity of $H(s)$, property \pref{G3} of the function $G(s)$ and the fact that the function $D(s)$ is bounded in the strip ${\mathcal S}_{\epsilon}$, we conclude that 
$H(s)$ is bounded on $\CC$ and $H(s)\to 0$ as $\im(s)\to \infty$.
Since $H$ is entire, it follows that $H \equiv 0$.

So far, we have proved that for all $s\in {\mathcal S}_{\epsilon}$ we have $f_1(s)h_2(s) = f_2(s) h_1(s)$. Let us define
$w(s):=f_1(s)/h_1(s)=f_2(s)/h_2(s)$. Since both $f_i(s)$ and $h_i(s)$ satisfy the same functional equation  
\eqref{matrix_equation}, if $s$ and $s+1$ belong to ${\mathcal S}_{\epsilon}$ we have
\begin{eqnarr*}
  w(s+1)h_1(s+1)
  &=& f_1(s+1) \\
  &=& A_{11}(s)f_1(s)+A_{12}(s)f_2(s) \\
  &=& w(s)[A_{11}(s)h_1(s)+A_{12}(s)h_2(s)],
\end{eqnarr*}
and therefore $w(s+1)=w(s)$.
Using again the fact that
$f_i$ and $h_i$ are  analytic in this strip and $h_i$ is also zero free there, 
we conclude that $w(s)$ is analytic in ${\mathcal S}_{\epsilon}$,
and the periodicity of $w$ implies that it may be extended to an entire periodic function. 
Lemma \ref{lemma1}(ii) together with the uniform boundedness of $f_i(s)$ in ${\mathcal S}_{\epsilon}$ imply 
 that there exists a constant $C>0$ such that for all $s\in {\mathcal S}_{\epsilon}$, 
\beqq
|w(s)|< C e^{\pi|\im(s)|}. 
\eeqq
By periodicity of $w$, we conclude that the above bound holds for all $s\in \c$. 
Since $w$ is periodic with period one, this bound implies that    
$w$ is a constant function (this follows from the Fourier series representation
of periodic analytic functions; see the proof of Proposition 2 in \cite{KP-HG}).
Finally, we know that $f_i(1)=h_i(1)=1$, and so we conclude that 
$w(s)\equiv 1$. Hence,  $f_i(s)\equiv h_i(s)$ for all $s\in {\mathcal S}_{\epsilon}$.
Since $h_i(s)$ are analytic in the wider strip\linebreak 
$-1/\alpha<\re(s)<2-1/\alpha$, by analytic continuation 
we conclude that \eqref{Mellin_transform} holds for all $s$ in $-1/\alpha<\re(s)<2-1/\alpha$.
\end{proof}

\begin{remark}
Since the proof of Theorem \ref{thm_Mellin_transform} is based on a verification technique,
it does not reveal how we derived the formula
on the right-hand side of \eqref{Mellin_transform}, for which we took a trial and error approach.
The expression
in \eqref{Mellin_transform} is already known, or may be easily computed, for
the spectrally positive case ($\rho=1-1/\alpha$; in this case
$T_0$ is the time of first passage below the level zero, and indeed
has a positive $1/\alpha$-stable law, as may be seen from \cite[Theorem 46.3]{Sato}),
the spectrally negative case ($\rho=1/\alpha$; due to \cite[Corollary 1]{Sim-hit})
and for the symmetric case ($\rho=1/2$; due to \cite{YYY-laws,Cor-thesis}),
we sought a function which interpolated these three cases and satisfied
the functional equation \eqref{matrix_equation}.
After a candidate was found, we verified that this was indeed the Mellin transform,
using the argument above.
\end{remark}

\newcommand*{\Lnorm}[1]{\lVert #1 \rVert}
\medskip\noindent
We turn our attention to computing the density of $T_0$.
Let us define $\Lnorm{x}=\min_{n \in \ZZ}\abs{x-n}$, and
\beqq
{\mathcal L}=\r \setminus (\q \cup \{x \in \r : \lim_{n\to \infty} \tfrac{1}{n}\ln\Lnorm{nx}=0\}). 
\eeqq
This set was introduced in \cite{ECP1601}, where it was shown that ${\mathcal L}$ is
a subset of the Liouville numbers 
and that
$x\in {\mathcal L}$ if and only if the coefficients
of the continued fraction representation of $x$ grow extremely fast.
It is known that ${\mathcal L}$ is dense, yet it is a rather small set:
it has Hausdorff dimension zero, and therefore its Lebesgue measure is also zero.

 For $\alpha \in \r$ we also define
 \beqq
 {\mathcal K}(\alpha)=\{N\in {\mathbb N} : \Lnorm{(N-\tfrac12)\alpha}>\exp(-\tfrac{\alpha-1}{2} (N-2)\ln(N-2)) \}.
 \eeqq
 \begin{proposition}\label{prop_K_alpha} Assume that $\alpha \notin \q$. 
 \begin{enumerate}[(i)]
 \item\label{i:prop_K_alpha:1} The set $ {\mathcal K}(\alpha)$ is unbounded and has density equal to one:
 \beqq
 \lim\limits_{n\to \infty} \frac{\# \{ {\mathcal K}(\alpha) \cap [1,n] \}}{n} =1.
 \eeqq 
 \item\label{i:prop_K_alpha:2} If $\alpha \notin {\mathcal L}$, the set  ${\mathbb N}\setminus {\mathcal K}(\alpha)$ is finite.  
 \end{enumerate}
 \end{proposition}
 \begin{proof}
Part \pref{i:prop_K_alpha:1} follows, after some short manipulation,
from the well-known fact that for any irrational
$\alpha$ the sequence $\Lnorm{(N-\tfrac12)\alpha}$ 
is uniformly distributed on the interval $(0,1/2)$.

To prove part \pref{i:prop_K_alpha:2}, first assume that $\alpha \notin {\mathcal L}$. 
Since $\lim_{n\to +\infty} \tfrac{1}{n}\ln\Lnorm{n \alpha}=0$, there exists $C>0$ such that 
for all $n$ we have $\Lnorm{n \alpha}>C2^{-n}$. Then for all $N$ we have
\beqq
\Lnorm{(N-\tfrac{1}{2})\alpha}\ge \frac{1}{2} \Lnorm{(2N-1)\alpha}>C2^{-2N}.
\eeqq
Since for all $N$ large enough it is true that 
\beqq
C2^{-2N}>\exp(-\tfrac{\alpha-1}{2} (N-2)\ln(N-2)),
\eeqq
we conclude that all $N$ large enough will be in the set ${\mathcal K}(\alpha)$, therefore 
the set ${\mathbb N}\setminus {\mathcal K}(\alpha)$ is finite. 
 \end{proof}

\begin{theorem}\label{t:T0 densities}
 Let $p$ be the density of $T_0$ under $\stP_1$.
 \begin{itemize}
 \item[(i)]
If $\alpha \notin \q$ then for all $t>0$ we have
\begin{IEEEeqnarray*}{rClll}
p(t)
  &=& \
  \lim_{\substack{N\in {\mathcal K}(\alpha)\\N\to \infty }}
  \Bigg[ &
  \frac{\sin\left(\frac{\pi}{\alpha} \right)}{\pi\sin(\pi \hat \rho)}
  \sum\limits_{1\le k < \alpha(N-\frac{1}{2})-1} &
  \sin(\pi \hat \rho(k+1))
  \frac{\sin\left(\frac{\pi}{\alpha} k \right)}{\sin\left( \frac{\pi}{\alpha}(k+1) \right)} \\
  &&&& {} \times \frac{\Gamma\left(\frac{k}{\alpha} +1\right) }{k! }
  (-1)^{k-1} t^{-1-\frac{k}{\alpha}} \\
  &&& \IEEEeqnarraymulticol{2}{l}{{} -
  \frac{\sin\left(\frac{\pi}{\alpha} \right)^2}{\pi\sin(\pi \hat \rho)}
  \sum\limits_{1\le k < N} 
  \frac{\sin(\pi \alpha \hat \rho k)}{\sin(\pi \alpha k)}
  \frac{\Gamma\left(k-\frac{1}{\alpha}\right)}{ \Gamma\left(\alpha k -1\right)}
  t^{-k-1+\frac{1}{\alpha}}\Bigg].}
\yesnumber \label{e:T0 density}
\end{IEEEeqnarray*}
The above limit is uniform for $t\in [\epsilon,\infty)$ and any $\epsilon>0$.  
\item[(ii)]
If $\alpha=m/n$ (where $m$ and $n$ are coprime natural numbers) then for all $t>0$ we have
\begin{eqnarr*}
\nonumber
p(t)
  &=&
  \frac{\sin\left(\frac{\pi}{\alpha} \right)}{\pi\sin(\pi \hat \rho)}
\sum\limits_{\substack{k\ge 1 \\ k\neq -1\; {\textnormal {(mod m)}}}} 
\sin(\pi \hat \rho(k+1))
\frac{\sin\left(\frac{\pi}{\alpha} k \right)}{\sin\left( \frac{\pi}{\alpha}(k+1) \right)}
\frac{\Gamma\left(\frac{k}{\alpha} +1\right) }{k! }
(-1)^{k-1} t^{-1-\frac{k}{\alpha}} \\ 
  && {} -
  \frac{\sin\left(\frac{\pi}{\alpha} \right)^2}{\pi\sin(\pi \hat \rho)}
\sum\limits_{\substack{k\ge 1 \\ k\neq 0\;{\textnormal {(mod n)}}}} 
\frac{\sin(\pi \alpha \hat \rho k)}{\sin(\pi \alpha k)}
\frac{\Gamma\left(k-\frac{1}{\alpha}\right)}{ \Gamma\left(\alpha k -1\right)}
t^{-k-1+\frac{1}{\alpha}} \\
  && {} - 
  \frac{\sin\left(\frac{\pi}{\alpha} \right)^2}{\pi^2 \alpha \sin(\pi \hat \rho)}
  \sum\limits_{k\ge 1}
  (-1)^{km} 
  \frac{\Gamma\left(kn-\frac{1}{\alpha}\right)}{(km-2)!}
  R_k (t)
  t^{-kn-1+\frac{1}{\alpha}}, \yesnumber
\label{e:T0 density2}
\end{eqnarr*}
where
\begin{eqnarr*}
  R_k (t) &:=&  \pi \alpha \hat \rho \cos(\pi \hat \rho km)\\
  &&  {} - \sin(\pi \hat \rho km)
  \left[\pi \cot\left(\tfrac{\pi}{\alpha}\right)
    - \psi\bigl(kn - \tfrac{1}{\alpha} \bigr)
    + \alpha \psi(km - 1)
    + \ln(t) \right]
\end{eqnarr*}
and $\psi$ is the digamma function.
    The three series in \eqref{e:T0 density2} converge uniformly for $t\in [\epsilon,\infty)$ and any $\epsilon>0$. 
\item[(iii)] 
For all values of $\alpha$ and any $c>0$, the following asymptotic expansion holds as $t\downto 0$:
\[
p(t) =
\frac{\alpha\sin\left(\frac{\pi}{\alpha} \right)}{\pi\sin(\pi \hat \rho)}
\sum\limits_{1\le n <1+c } 
\sin(\pi \alpha \hat \rho n) 
\frac{\Gamma(\alpha n+1)}{ \Gamma\left(n+\frac{1}{\alpha} \right)}
(-1)^{n-1} t^{n-1+\frac{1}{\alpha}}+
O(t^{c+\frac{1}{\alpha}}).
\]
\end{itemize}
\begin{proof}
Recall that $h_1(s)=\stE_1[T_0^{s-1}]$ denotes the function in \eqref{Mellin_transform}. According to 
Lemma \ref{lemma1}(ii), for any $x\in \r$, $h_1(x+\iu y)$ decreases to zero exponentially
fast as $y\to \infty$. This implies that the density of $T_0$ exists and is a smooth function. It also implies that $p(t)$ can be written as the inverse Mellin transform,
\beq\label{px_inverse_Mellin}
p(t)=\frac{1}{2\pi \iu} \int\limits_{1+\iu \r} h_1(s) t^{-s} \d s.  
\eeq
The function $h_1(s)$ is meromorphic, and it has poles at points
\beqq
\{s^{(1)}_n:=1+n/\alpha\; : \; n\ge 1\} \cup \{s^{(2)}_n:=n-1/\alpha \; : \; n\ge 2\}
\cup\{s^{(3)}_n:=-n-1/\alpha \;:\; n\ge 0\}
\eeqq
If $\alpha \notin \q$, all these points are distinct and $h_1(s)$ has only simple poles. When $\alpha \in \q$, 
some of $s^{(1)}_n$ and $s^{(2)}_m$ will coincide, and $h_1(s)$ will have double poles at these points. 

Let us first consider the case $\alpha \notin \q$, so that all poles are simple. 
Let $N\in {\mathcal K}(\alpha)$ and define $c=c(N)=N+\tfrac{1}{2}-\tfrac{1}{\alpha}$. 
Lemma \ref{lemma1}(ii) tells us that $h_1(s)$ decreases exponentially to zero as $\im(s)\to \infty$,
uniformly in any finite vertical strip. Therefore, we can shift the contour of integration in  
\eqref{px_inverse_Mellin} and obtain, by Cauchy's residue theorem,
\begin{eqnarr}
\nonumber
p(t)
  &=&
  - \sum_n
  \Res_{s=s^{(1)}_n}(h_1(s)t^{-s})
  - \sum_m
  \Res_{s=s^{(2)}_m}(h_1(s)t^{-s}) \\
  && {} + \frac{1}{2\pi \iu} \int\limits_{c(N)+\iu \r} h_1(s) t^{-s} \d s,
\label{px_inverse_Mellin2}
\end{eqnarr}
where $\sum_n$ and $\sum_m$ indicate
summation over $n \ge 1$ such that $s_n^{(1)} < c(N)$
and over $m \ge 2$ such that $s_m^{(2)} < c(N)$, respectively.
Computing the residues we obtain the two sums in the right-hand side of  \eqref{e:T0 density}. 

Now our goal is to show that the integral term
\beqq
I_N(t):=\frac{1}{2\pi \iu} \int\limits_{c(N)+\iu \r} h_1(s) t^{-s} \d s
\eeqq 
converges to zero as $N\to +\infty$, $N\in {\mathcal K}(\alpha)$. We use the reflection formula 
for the gamma function and the inequalities
\begin{IEEEeqnarray*}{l"l}
  \abs{\sin(\pi x)} > \Lnorm{x}, & x \in \RR, \\
  \abs{\sin(x)}\cosh(y) \le \abs{\sin(x+\iu y)} = \sqrt{\cosh^2(y)-\cos^2(x)} \le \cosh(y),
  & x,y\in \r,
\end{IEEEeqnarray*}
to estimate $h_1(s)$, $s=c(N)+\iu u$, as follows 
\beq\label{estimate1}
\abs{h_1(s)}&=&
\frac{\sin\left(\frac{\pi}{\alpha} \right)}{\sin(\pi \hat \rho)}
\biggabs{ \frac{\sin\left( \pi \hat \rho (1-\alpha+\alpha s)\right)}
    {\sin\left(\frac{\pi}{\alpha} (1-\alpha+\alpha s) \right)}
    \frac{\sin(\pi s)}{\sin(\pi  \alpha (s-1))}
  \frac{\Gamma(s-1)}{\Gamma(\alpha (s-1))} } \nonumber \\ 
  &  \le&  \frac{C_1}{\Lnorm{\alpha (N-\frac{1}{2})}}
\frac{  \cosh(\pi \alpha \hat \rho u) }{\cosh(\pi \alpha u)}
\biggabs{
  \frac{\Gamma(s-1)}{\Gamma(\alpha (s-1))} } .
\eeq

Using Stirling's formula (\ref{Stirlings_formula}),
we find that 
\beq\label{e:ratio of gammas}
\frac{\Gamma(s)}{\Gamma(\alpha s)}=\sqrt{\alpha} 
e^{-s \left((\alpha-1)\ln(s)+A\right)+O(s^{-1})}, \for s\to \infty, \; \re(s)>0,
\eeq
where $A:=1-\alpha+\alpha\ln(\alpha)>0$.
Therefore, there exists a constant $C_2>0$ such that for $\re s >0$  we can estimate
\begin{eqnarr*}
\biggabs{ \frac{\Gamma(s)}{\Gamma(\alpha s)}}
  &<& C_2 e^{- (\alpha-1) \re (s) \ln(\Re(s)) +(\alpha-1) \abs{\Im(s)} \frac{\pi}{2} }.
\end{eqnarr*}
Combining the above estimate with \eqref{estimate1} and using the fact that $N\in {\mathcal K}(\alpha)$ we find that 
\begin{eqnarr*}
  \abs{h_1(c(N)+\iu u)}
  &<& \frac{C_1C_2}{\Lnorm{\alpha (N-\frac{1}{2})}}
  \frac{  \cosh(\pi \alpha \hat \rho u) }{\cosh(\pi \alpha u)}
  e^{- (\alpha-1) (c(N)-1) \ln(c(N)-1) +(\alpha-1) \abs{u} \frac{\pi}{2} } \\
  &<& C_1C_2 e^{-\frac{\alpha-1}{2} (N-2)\ln(N-2)} 
  \frac{  \cosh(\pi \alpha \hat \rho u) }{\cosh(\pi \alpha u)} e^{(\alpha-1) \abs{u} \frac{\pi}{2} } .
\end{eqnarr*}
Note that the function in the right-hand side of the above inequality decreases to zero exponentially fast 
as $\abs{u}\to \infty$ (since $\alpha \hat \rho + \tfrac{1}{2} (\alpha-1)-\alpha<0$),
and hence in particular is integrable on $\RR$.
Thus we can estimate 
\begin{multline*}
\abs{I_N(t)}
  = \frac{t^{-c(N)}}{2\pi}
  \left\lvert\int_\r h_1(c(N)+\iu u)) t^{-\iu u} \, \d u \right\rvert
  < \frac{t^{-c(N)}}{2\pi} \int_\r |h_1(c(N)+\iu u))| \, \d u \\
  {} < \frac{t^{-c(N)}}{2\pi}C_1C_2 e^{-\frac{\alpha-1}{2} (N-2)\ln(N-2)} 
\int_\r \frac{  \cosh(\pi \alpha \hat \rho u) }{\cosh(\pi \alpha u)} e^{(\alpha-1) \abs{u} \frac{\pi}{2} } \, \d u .
\end{multline*}
When $N\to \infty$, the quantity in the right-hand side of the above inequality
converges to zero for every $t>0$. This ends the proof of part (i). 

The proof of part (ii) is very similar, and we offer only a sketch.
It also begins with \eqref{px_inverse_Mellin2}
and uses the above estimate for $h_1(s)$. The only difference is that when $\alpha \in \q$
some of $s^{(1)}_n$ and $s^{(2)}_m$ will coincide, and $h_1(s)$ will have
double poles at these points.
The terms with double poles give rise to the third series in \eqref{e:T0 density2}. In this case all
three series are convergent, and we can express the limit of partial sums as
a series in the usual sense.

The proof of part (iii) is much simpler: we need to shift the contour of integration in 
\eqref{px_inverse_Mellin2} in the opposite direction ($c<0$). The proof is identical to the
proof of Theorem 9 in \cite{Kuz2011}. 
\end{proof}
\end{theorem}

\begin{remark}
  We offer some remarks on the asymptotic expansion. 
  When $\rhohat = 1/\alpha$, all of its terms are equal to zero. This is the spectrally
  positive case, in which $T_0$ has the law of a positive $1/\alpha$-stable
  random variable, and it is known
  that its density is exponentially small at zero; see \cite[equation (14.35)]{Sato}
  for a more precise result.
  
  Otherwise, the series given by including all the terms in (iii) is divergent for all $t > 0$.
  This may be seen from the fact that the terms do not approach $0$;
  we sketch this now. When $\alpha\rhoh \in \QQ$, some terms are zero,
  and in the rest the sine term is bounded away from zero;
  when $\alpha\rhoh \notin \QQ$, it follows that
  $\limsup_{n \to \infty} \abs{\sin(\pi\alpha\rhoh n)} = 1$.
  One then bounds the
  ratio of gamma functions from below by $\Gamma(\alpha n)/\Gamma((1+\epsilon)n)$,
  for some small enough $\epsilon > 0$ and large $n$.
  This grows superexponentially due to \eqref{e:ratio of gammas},
  so $t^n\Gamma(\alpha n+1)/\Gamma(n+1/\alpha)$ is unbounded as $n \to \infty$.
\end{remark}

The next corollary shows that, for almost all irrational $\alpha$, the expression
\eqref{e:T0 density} can be written in a simpler form. 
\begin{corollary}\label{cor_alpha_notin_L}
If $\alpha \notin {\mathcal L} \cup \q$ then 
\begin{eqnarr}
p(t)&=&\frac{\sin\left(\frac{\pi}{\alpha} \right)}{\pi\sin(\pi \hat \rho)}
\sum\limits_{k\ge 1}
\sin(\pi \hat \rho(k+1))
\frac{\sin\left(\frac{\pi}{\alpha} k \right)}{\sin\left( \frac{\pi}{\alpha}(k+1) \right)}
\frac{\Gamma\left(\frac{k}{\alpha} +1\right) }{k! }
(-1)^{k-1} t^{-1-\frac{k}{\alpha}} \nonumber \\
&& {} -
\frac{\sin\left(\frac{\pi}{\alpha} \right)^2}{\pi\sin(\pi \hat \rho)}
\sum\limits_{k\ge 1}
\frac{\sin(\pi \alpha \hat \rho k)}{\sin(\pi \alpha k)}
\frac{\Gamma\left(k-\frac{1}{\alpha}\right)}{ \Gamma\left(\alpha k -1\right)}
t^{-k-1+\frac{1}{\alpha}}.
\label{e:T0 density4}
\end{eqnarr}
The two series in the right-hand side of the above formula converge uniformly for
$t\in [\epsilon,\infty)$ and any $\epsilon>0$. 
\end{corollary}
\begin{proof}
As we have shown in Proposition \ref{prop_K_alpha}, if $\alpha \notin {\mathcal L} \cup \q$ then 
the set ${\mathbb N} \setminus{\mathcal K}(\alpha)$ is finite. Therefore we can remove the restriction 
$N \in {\mathcal K}(\alpha)$ in \eqref{e:T0 density}, and need only show that
both series in \eqref{e:T0 density4} converge. 


In \cite[Proposition 1]{ECP1601} it was shown that $x \in {\mathcal L}$ iff $x^{-1} \in {\mathcal L}$. Therefore,
according to our assumption, 
both $\alpha$ and $1/\alpha$ are not in the set ${\mathcal L}$.
From the definition of ${\mathcal L}$ we see that there exists $C>0$ such that
 $\Lnorm{\alpha n}>C2^{-n}$ and $\Lnorm{\alpha^{-1} n}>C2^{-n}$ for all integers $n$. 
 Using the estimate $|\sin(\pi x)|\ge \Lnorm{x}$ and 
 Stirling's formula \eqref{Stirlings_formula}, it is easy to see that both series in 
 in  \eqref{e:T0 density4} converge (uniformly for $t\in [\epsilon,\infty)$ and any $\epsilon>0$),
 which ends the proof of the corollary.  
\end{proof}

\begin{remark}
Note that formula \eqref{e:T0 density4} may not be true if $\alpha \in {\mathcal L}$, as the series 
may fail to converge. An example where this occurs is given after Theorem 5 in \cite{Kuz2011}.
\end{remark}

\section{Applications}
\label{s:app}

\subsection{Conditioning to avoid zero}

In \cite[\S 4.2]{CPR}, Chaumont, Pant{\'{\i}} and Rivero
discuss a harmonic transform
of a stable process with $\alpha > 1$
which results in \define{conditioning to avoid zero}. The
results quoted in that paper are a special case of the
notion of
conditioning a L\'evy process to avoid zero, which is explored
in \citet{Pan-cond}.

In these works, in terms of the parameters used in the introduction, the authors define
\begin{equation}
\label{e:h}
  h(x) = 
  \frac{\Gamma(2-\alpha) \sin(\pi\alpha/2)}
  {c \pi (\alpha-1) (1+\beta^2 \tan^2(\pi\alpha/2))}
   (1-\beta \sgn(x)) \abs{x}^{\alpha-1},
  \for x \in \RR.
\end{equation}
If we write the function $h$ in terms of the $(\alpha,\rho)$
parameterisation which we prefer, this gives
\[ h(x) =
  - \Gamma(1-\alpha)
  \frac{\sin (\pi\alpha\rhoh)}{\pi}
  \abs{x}^{\alpha-1} , \for x > 0, \]
and the same expression with $\rhoh$ replaced by $\rho$ when $x < 0$.

In \cite{Pan-cond}, Pant{\'\i} proves the following 
proposition for all L\'evy processes, and $x \in \RR$,
with a suitable
definition of $h$. Here we quote only the result for stable
processes and $x \ne 0$.
Hereafter, $\FFt$ is the standard filtration associated with $X$,
and
$n$ refers to the excursion measure of the stable process
away from zero, normalised
(see \cite[(7)]{Pan-cond} and \cite[(4.11)]{FG-bridges}) such that
\[ n(1-e^{-q\zeta}) = 1/u^q(0), \]
where $\zeta$ is the excursion length and $u^q$ is the $q$-potential
density of the stable process.
\def\tempPprop{\cite[Theorem 2, Theorem 6]{Pan-cond}}
\begin{proposition}[\tempPprop]\label{p:Panti}
Let $X$ be a stable process, and $h$ the function in \eqref{e:h}.
\begin{enumerate}[(i)]
\item \label{p:Panti:1}
The function $h$ is invariant for the stable process
killed on hitting $0$, that is,
\begin{equation}\label{e:h invariant}
  \stE_x[ h(X_t), t < T_0] = h(x), \for t > 0, \, x \neq 0.
\end{equation}
Therefore, we may define a family of measures $\stPaz_x$ by
\[ \stPaz_x(\Lambda)
  = \frac{1}{h(x)} \stE_x[ h(X_t) \Ind_\Lambda, t < T_0 ],
  \for x \ne 0, \, \Lambda \in \FF_t, \]
for any $t \ge 0$.

\item \label{p:Panti:2}
The function $h$ can be represented as
\begin{equation*}
  h(x) = \lim_{q \downto 0} \frac{\stP_x(T_0 > \ee_q)}{n(\zeta > \ee_q)},
  \for x \ne 0,
\end{equation*}
where $\ee_q$ is an independent exponentially distributed random
variable with parameter $q$. Furthermore, for any stopping time $T$
and $\Lambda \in \FF_T$, and any $x \ne 0$,
\[ \lim_{q \downto 0} \stP_x(\Lambda, T < \ee_q \vert T_0 > \ee_q) = \stPaz_x(\Lambda) . \]
\end{enumerate}
\end{proposition}
This justifies the name `the stable process conditioned to avoid zero'
for the canonical process associated
with the measures $(\stPaz_x)_{x \ne 0}$.
We will denote this process by $\Xaz$.

\medskip\noindent
Our aim in this section is to prove the following variation
of Proposition
\ref{p:Panti}(ii), making use of our expression for the density of $T_0$.
Our presentation here owes much to
\citet[\S 4.3]{YYY-pen}.

\begin{proposition}\label{p:lim P} Let $X$ be a stable process adapted to the filtration $\FFt$,
and $h \colon \RR \to \RR$ as in \eqref{e:h}. 
\begin{enumerate}[(i)]
\item \label{p:lim P:1}
Define the function
\begin{equation*} Y(s,x) = \frac{\stP_x(T_0 > s)}{h(x) n(\zeta > s)}
  \for s > 0, \, x \ne 0. 
  \end{equation*}
Then, for any $x \ne 0$,
\begin{equation}
\lim_{s \to \infty} Y(s,x) = 1,
  \label{YYY}
\end{equation}
and furthermore, $Y$ is bounded away from $0$ and $\infty$ on its whole domain.

\item \label{p:lim P:2}
For any $x \ne 0$,
stopping time $T$ such that $\LevE_x[T] < \infty$, and
$\Lambda \in \FF_T$,
\[ \stPaz_x(\Lambda)
  = \lim_{s \to \infty} \stP_x (\Lambda \vert T_0 > T + s) . \]
\end{enumerate}

\begin{proof}
We begin by proving
\begin{equation}\label{e:h lim s}
h(x) = \lim_{s \to \infty} \frac{\stP_x(T_0 > s)}{n(\zeta > s)},
\end{equation}
for $x > 0$, noting that when $x < 0$, we may deduce the same limit by duality.

Let us denote the density of the measure $\stP_x(T_0 \in \cdot)$
by $p(x,\cdot)$. A straightforward application of scaling
shows that
\[ \stP_x(T_0 > t) = \stP_1(T_0 > x^{-\alpha} t), \for x > 0, \, t \ge 0, \]
and so we may focus our attention on $p(1,t)$, which
is the quantity given as $p(t)$ in Theorem \ref{t:T0 densities}.
In particular, we have
\[
  p(1,t) = - \frac{ \sin^2 (\pi/\alpha)}{\pi \sin(\pi \rhoh)}
  \frac{\sin(\pi\alpha\rhoh)}{\sin(\pi\alpha)}
  \frac{\Gamma(1-1/\alpha)}{\Gamma(\alpha-1)}
  t^{1/\alpha-2}
  + O(t^{-1/\alpha-1})
  .
\]
Denote the coefficient of $t^{1/\alpha-2}$ in the
first term of this expression by $P$.


To obtain an expression for $n(\zeta > t)$, we turn to
\citet{FG-bridges}, in which the authors compute
explicitly the density
of $n(\zeta \in \cdot)$ for a stable process; see p.~84
in that work, where $n$ is denoted $P^*$ and $\zeta$
is denoted $R$.
The authors work with
a different normalisation of the stable process; they have
$c = 1$. In our context, their result says
\begin{equation}\label{e:n density}
  n(\zeta \in \dd t)
  = \frac{\alpha-1}{\Gamma(1/\alpha)}
  \frac{\sin(\pi/\alpha)}{\cos(\pi(\rho-1/2))}
  t^{1/\alpha-2} \, \dd t, \qquad t\geq 0.
\end{equation}
Denote the coefficient in the above power law by $W$.

We can now compute $h$. We will  use elementary
properties of
trigonometric functions and the reflection identity for the gamma
function. For $x > 0$,
\begin{eqnarr*}
  \frac{\stP_x(T_0 > t)}{n(\zeta > t)}
  &=& \frac{P}{W} x^{\alpha-1}
  + O(t^{1-2/\alpha}) \\
  &=&
  - \frac{\cos(\pi(\rho-1/2)) \sin(\pi\alpha\rhoh)}
  {\Gamma(\alpha) \sin(\pi\rhoh) \sin(\pi\alpha)} x^{\alpha-1}
  + o(1) \\
  &=&
  - \frac{1}{\Gamma(\alpha)}
  \frac{\sin(\pi\alpha\rhoh)}{\sin(\pi\alpha)} x^{\alpha-1}
  + o(1) \\
  &=& - \Gamma(1-\alpha) \frac{\sin(\pi\alpha\rhoh)}{\pi} x^{\alpha-1} 
  + o(1).
\end{eqnarr*}
This proves \eqref{e:h lim s} for $x > 0$, and it is simple to deduce
via duality that this limit holds for $x \ne 0$.

We now turn our attention to the slightly more delicate result about $Y$.
It is clear that the limit in \eqref{YYY} holds, so we only need to prove that $Y$ is bounded.
We begin by noting that, for fixed $x \ne 0$, $Y(t,x)$ is bounded 
away from $0$ and $\infty$ in $t$
since the function is continuous and converges to $1$ as $t \to \infty$.
Now, due to the expression \eqref{e:n density} and the scaling property of
$X$, we have the relation $Y(t,x) = Y(\abs{x}^{-\alpha}t,\sgn x)$.
This then shows that $Y$ is bounded as a function of two variables.

With this in hand, we move on to the calculation of the limiting measure.
This proceeds along familiar lines, using the strong Markov property:
\begin{eqnarr*}
  \stP_x( \Lambda \vert T_0 > T+s)
  &=& \stE_x \biggl[ \frac{ \stP_x(\Ind_\Lambda, T_0 > T+s \vert \FF_T) }
    {\stP_x(T_0 > T+s)} \biggr] \\
  &=& \stE_x\biggl[ \Ind_\Lambda \Indic{T_0 > T}
    \frac{\stP_{X_T}(T_0 > s)}{\stP_x(T_0 > T+s)} \biggr] \\
  &=& \stE_x\biggl[ \Ind_\Lambda \Indic{T_0 > T} h(X_T)
  Y(s, X_T)
  \frac{n(\zeta > s)}{n(\zeta > T+s)}
  \frac{1}{h(x) Y(s+T,x)}
  \biggr].
\end{eqnarr*}
Now, as $h$ is invariant for the stable process killed at zero,
\eqref{e:h invariant} also holds at $T$, and in particular
the random variable $h(X_T)\Indic{T_0 > T}$ is integrable; meanwhile, $Y$ is bounded
away from zero and $\infty$.
\newcommand\notfrac[2]{#1/#2}
Finally, $\notfrac{n(\zeta > s)}{n(\zeta>T+s)} = \bigl(1+\frac{T}{s}\bigr)^{1-1/\alpha}$;
the moment condition in the statement of the theorem permits us to
apply the dominated convergence theorem, and this gives the result.
\end{proof}
\end{proposition}

\medskip \noindent
We offer a brief comparison to conditioning a L\'evy process to stay
positive. In this case, \citet[Remark 1]{Cha-cond} observes that
the analogue of
Proposition \refpref{p:lim P}{p:lim P:2}
holds under Spitzer's condition, and in particular for a stable process.
However, it appears that in general, a key role is played by
the exponential random variable
analogous to that
appearing in Proposition \refpref{p:Panti}{p:Panti:2}.

\subsection{The radial part of the stable process conditioned to avoid zero}
\label{ss:Raz}

For this section, consider $X$ to be the symmetric stable process,
as in section \ref{s:symm}. There we computed the Lamperti transform $\xi$
of the pssMp
\[ R_t = \abs{X_t}\Indic{t < T_0}, \for t \ge 0, \]
and gave its characteristic exponent $\Psi$ in \eqref{e:WHF low index}.

Consider now the process
\[ \Raz_t = \bigabs{\Xaz_t}, \for t \ge 0.\]
This is also
a pssMp, and we may consider its Lamperti transform, which we will
denote by $\xiaz$. The characteristics of the Lamperti--Kiu
representation of $\Xaz$ have been computed explicitly in \cite{CPR},
and the characteristic exponent, $\CEaz$, of $\xiaz$ could be
computed from this information; however, the harmonic transform
in Proposition \refpref{p:Panti}{p:Panti:1}
gives us the following straightforward relationship
between characteristic exponents:
\[ \CEaz(\theta) = \CE(\theta - \iu(\alpha-1)) .  \]
This allows us to calculate
\[ \CE^\updownarrow(\theta)
  = 2^\alpha \frac{\Gamma(1/2-\iu\theta/2)}{\Gamma((1-\alpha)/2-\iu\theta/2)}
  \frac{\Gamma(\alpha/2 + \iu\theta/2)}{\Gamma(\iu\theta/2)} , \for \theta \in \RR.
\]
It is immediately apparent that $\xiaz$ is the dual L\'evy process
to $\xi$. It then follows that $R$ is a time-reversal
of $\Raz$, in the sense of \cite[\S 2]{CP-lower}: roughly speaking,
if one fixes $x > 0$,
starts the process $\Raz$ at zero (as in \cite{CC-wc}, say)
and runs it backward from its last moment below some level $y$, where $y>x$,
simultaneously conditioning on the position of the left limit at this time taking value $x$,
then one obtains the law of $R$ under $\stP_x$.

We remark that this relationship is already known for Brownian
motion, where $\Raz$ is a Bessel process of dimension $3$.
However, it seems unlikely that any such time-reversal
property will hold for a general L\'evy process conditioned to avoid
zero.

\subsection{The entrance law of the excursion measure}

It is known, from a more general result \cite[(2.8)]{CFY-ext}
on Markov processes in weak duality, that for any Borel function
$f$, the equality
\[
  \int_0^\infty e^{-qt} n(f(X_t)) \, \dd t
  = \int_{\RR} f(x) \stEh_x\bigl[ e^{-qT_0} \bigr] \, \dd x \]
holds, where $n$ is the excursion measure of $X$ from zero.
(This formulation is from \cite[(3.9)]{YYY-pen}.) In terms of
densities, this may be written
\begin{eqnarr*}
  n(X_t \in \dd x) \, \dd t
  &=& \stPh_x(T_0 \in \dd t) \, \dd x \\
  &=& \abs{x}^{-\alpha} \stP_{\sgn(-x)}(T_0 \in \abs{x}^{-\alpha}\dd t) \, \dd x
\end{eqnarr*}
Therefore, our expressions in Theorem \ref{t:T0 densities}
for the density of $T_0$ yield expressions
for the density of the entrance law of the excursions of the stable process
from zero.

\subsection{Identities in law using the exponential functional}

In a series of papers (\citet{BY-ent}, \citet{CC-wc}, \citet{CKPR-fluc})
it is proved that under certain conditions, the laws $(\stP_x)_{x > 0}$ of an
$\alpha$-pssMp $X$ admit a weak limit $\stP_0$ as $x \downto 0$, in the Skorokhod space
of c\`adl\`ag paths.
If $\xi$ is the Lamperti
transform of $X$ under $\stP_1$,
then provided that $\LevE\abs{\xi_1} < \infty$
and $m := \LevE \xi_1 > 0$, it is known that the entrance law of $\stP_0$
satisfies
\[ \stE_0( f(X_t^\alpha)) = \frac{1}{\alpha m} \LevE( I(-\alpha\xi)^{-1} f(t/I(-\alpha\xi))) , \]
for any $t > 0$ and Borel function $f$.
Similar expressions are available under
less restrictive conditions on $\xi$.

It is tempting to speculate that any rssMp may admit a weak limit
$\stP_0$ along similar lines, but we do not propose any results
in this direction; instead, we demonstrate
similar formulae for the entrance law $n(X_t \in \cdot)$ of the
stable process, and the corresponding measure $\stPaz_0$ for the stable
process conditioned to avoid zero.

\medskip \noindent
Let $X$ be a stable process, possibly asymmetric.
From the previous subsection, we have that
\[ n(f(X_t))
  = \int_{-\infty}^\infty
  \abs{x}^{-\alpha} p(\sgn(-x), \abs{x}^{-\alpha}t) f(x) \, \dd x.
\]
Substituting in the integral, and recalling that the law of $T_0$
for the stable process is equal to the law of the exponential
functional $I(\alpha\xi)$ of the Markov additive process associated with it,
we obtain
\begin{eqnarr*}
  n(f(X_t))
  &=& \frac{1}{\alpha}\int_0^\infty p(1,u)
  f( -(u/t)^{-1/\alpha}) u^{-1/\alpha} t^{1/\alpha-1} \, \dd u \\
  && {} + \frac{1}{\alpha}\int_0^\infty p(-1,u)
  f( (u/t)^{-1/\alpha}) u^{-1/\alpha} t^{1/\alpha-1} \, \dd u \\
  &=& \frac{1}{\alpha}\MAPE_1\bigl[ f(-(t/I(\alpha\xi))^{1/\alpha})
  I(\alpha\xi)^{-1/\alpha} t^{1/\alpha-1} \bigr] \\
&&  {} + \frac{1}{\alpha}\MAPE_2\bigl[ f((t/I(\alpha\xi))^{1/\alpha})
  I(\alpha\xi)^{-1/\alpha} t^{1/\alpha-1} \bigr].
\end{eqnarr*}

Recall from \cite{Pan-cond} that the
law $\stPaz_0$ of the stable process conditioned to avoid
zero is given by the following harmonic transform of the
stable excursion measure $n$:
\[ \stPaz_0(\Lambda) = n( \Ind_\Lambda h(X_t), t < \zeta),
  \for t \ge 0, \, \Lambda \in \FF_t , \]
with $h$ as in \eqref{e:h}. Therefore, applying the above result to the Borel function $hf$, we
obtain
\begin{eqnarr*}
\stEaz_0(f(X_t))
  &=& n( h(X_t) f(X_t) ) \\
  &=& 
  \Gamma(-\alpha)\frac{\sin(\pi\alpha\rho)}{\pi}
  \MAPE_1\bigl[I(\alpha\xi)^{-1} f(-(t/I(\alpha\xi))^{1/\alpha})
  \bigr] \\
  &&{} + \Gamma(-\alpha)\frac{\sin(\pi\alpha\rhoh)}{\pi}
  \MAPE_2\bigl[ I(\alpha\xi)^{-1} f((t/I(\alpha\xi))^{1/\alpha})
  \bigr] ,
\end{eqnarr*}
where we emphasise that $I(\alpha\xi)$ (under $\MAPE_i$)
is the exponential functional of the
Markov additive process associated to $X$.


\bigskip
\ACKNO{%
  Some of this research was carried out whilst AEK and ARW were visiting
  ETH Z\"urich and CIMAT. Both authors would like to offer thanks to both
  institutions for their hospitality. AK acknowledges the support by the
  Natural Sciences and Engineering Research Council of Canada. JCP acknowledges
  the support by  CONACYT (grant 128896). 
  All authors would like to thank the referee for their detailed comments that led
  to an improved version of this paper.%
}


\begin{thebibliography}{39}
\providecommand{\natexlab}[1]{#1}
\providecommand{\url}[1]{\texttt{#1}}
\expandafter\ifx\csname urlstyle\endcsname\relax
  \providecommand{\doi}[1]{doi: #1}\else
  \providecommand{\doi}{doi: \begingroup \urlstyle{rm}\Url}\fi
\providecommand{\MR}[1]{MR#1}
\providecommand{\arxivref}[2]{#1\ifx&#2&\else\ [#2]\fi}

\bibitem[Andrews et~al.(1999)Andrews, Askey, and Roy]{AAR}
G.~E. Andrews, R.~Askey, and R.~Roy.
\newblock \emph{Special functions}, volume~71 of \emph{Encyclopedia of
  Mathematics and its Applications}.
\newblock Cambridge University Press, Cambridge, 1999.
\newblock ISBN 0-521-62321-9; 0-521-78988-5. \MR{1688958 (2000g:33001)}

\bibitem[Artin(1964)]{Art-gamma}
E.~Artin.
\newblock \emph{The gamma function}.
\newblock Translated by Michael Butler. Athena Series: Selected Topics in
  Mathematics. Holt, Rinehart and Winston, New York, 1964. \MR{0165148 (29
  \#2437)}

\bibitem[Asmussen(2000)]{Asm-rp1}
S.~Asmussen.
\newblock \emph{Ruin probabilities}, volume~2 of \emph{Advanced Series on
  Statistical Science \& Applied Probability}.
\newblock World Scientific Publishing Co. Inc., River Edge, NJ, 2000.
\newblock ISBN 981-02-2293-9. \MR{1794582 (2001m:62119)}

\bibitem[Asmussen(2003)]{Asm-apq2}
S.~Asmussen.
\newblock \emph{Applied probability and queues}, volume~51 of
  \emph{Applications of Mathematics (New York)}.
\newblock Springer-Verlag, New York, second edition, 2003.
\newblock ISBN 0-387-00211-1. \MR{1978607 (2004f:60001)}

\bibitem[Ben-Ari and Neumann(2012)]{BN-Perron}
I.~Ben-Ari and M.~Neumann.
\newblock Probabilistic approach to {P}erron root, the group inverse, and
  applications.
\newblock \emph{Linear and Multilinear Algebra}, 60\penalty0 (1):\penalty0
  39--63, 2012.
\newblock ISSN 0308-1087. \MR{2869672}

\bibitem[Bertoin(1996)]{BertoinLP}
J.~Bertoin.
\newblock \emph{L\'evy processes}, volume 121 of \emph{Cambridge Tracts in
  Mathematics}.
\newblock Cambridge University Press, Cambridge, 1996.
\newblock ISBN 0-521-56243-0. \MR{1406564 (98e:60117)}

\bibitem[Bertoin and Yor(2002)]{BY-ent}
J.~Bertoin and M.~Yor.
\newblock The entrance laws of self-similar {M}arkov processes and exponential
  functionals of {L}\'evy processes.
\newblock \emph{Potential Anal.}, 17\penalty0 (4):\penalty0 389--400, 2002.
\newblock ISSN 0926-2601. \MR{1918243 (2003i:60082)}

\bibitem[Blumenthal and Getoor(1968)]{BG-mppt}
R.~M. Blumenthal and R.~K. Getoor.
\newblock \emph{Markov processes and potential theory}.
\newblock Pure and Applied Mathematics, Vol. 29. Academic Press, New York,
  1968. \MR{0264757 (41 \#9348)}

\bibitem[Caballero and Chaumont(2006{\natexlab{a}})]{CC-wc}
M.~E. Caballero and L.~Chaumont.
\newblock Weak convergence of positive self-similar {M}arkov processes and
  overshoots of {L}\'evy processes.
\newblock \emph{Ann. Probab.}, 34\penalty0 (3):\penalty0 1012--1034,
  2006{\natexlab{a}}.
\newblock ISSN 0091-1798. \MR{2243877 (2008c:60035)}

\bibitem[Caballero and Chaumont(2006{\natexlab{b}})]{CC06}
M.~E. Caballero and L.~Chaumont.
\newblock Conditioned stable {L}\'evy processes and the {L}amperti
  representation.
\newblock \emph{J. Appl. Probab.}, 43\penalty0 (4):\penalty0 967--983,
  2006{\natexlab{b}}.
\newblock ISSN 0021-9002. \MR{2274630 (2008d:60058)}

\bibitem[Caballero et~al.(2011)Caballero, Pardo, and P{\'e}rez]{CPP-radial}
M.~E. Caballero, J.~C. Pardo, and J.~L. P{\'e}rez.
\newblock Explicit identities for {L}\'evy processes associated to symmetric
  stable processes.
\newblock \emph{Bernoulli}, 17\penalty0 (1):\penalty0 34--59, 2011.
\newblock ISSN 1350-7265. \MR{2797981 (2012b:60186)}

\bibitem[Carmona et~al.(1997)Carmona, Petit, and Yor]{CPY97}
P.~Carmona, F.~Petit, and M.~Yor.
\newblock On the distribution and asymptotic results for exponential
  functionals of {L}\'evy processes.
\newblock In \emph{Exponential functionals and principal values related to
  {B}rownian motion}, Bibl. Rev. Mat. Iberoamericana, pages 73--130. Rev. Mat.
  Iberoamericana, Madrid, 1997. \MR{1648657 (99h:60144)}

\bibitem[Chaumont(1996)]{Cha-cond}
L.~Chaumont.
\newblock Conditionings and path decompositions for {L}\'evy processes.
\newblock \emph{Stochastic Process. Appl.}, 64\penalty0 (1):\penalty0 39--54,
  1996.
\newblock ISSN 0304-4149. \MR{1419491 (98b:60131)}

\bibitem[Chaumont and Pardo(2006)]{CP-lower}
L.~Chaumont and J.~C. Pardo.
\newblock The lower envelope of positive self-similar {M}arkov processes.
\newblock \emph{Electron. J. Probab.}, 11:\penalty0 no. 49, 1321--1341, 2006.
\newblock ISSN 1083-6489. \MR{2268546 (2008f:60042)}

\bibitem[Chaumont et~al.(2012)Chaumont, Kyprianou, Pardo, and
  Rivero]{CKPR-fluc}
L.~Chaumont, A.~Kyprianou, J.~C. Pardo, and V.~Rivero.
\newblock Fluctuation theory and exit systems for positive self-similar
  {M}arkov processes.
\newblock \emph{Ann. Probab.}, 40\penalty0 (1):\penalty0 245--279, 2012.
\newblock ISSN 0091-1798. \MR{2917773}

\bibitem[Chaumont et~al.(2013)Chaumont, Pant{\'{\i}}, and Rivero]{CPR}
L.~Chaumont, H.~Pant{\'{\i}}, and V.~Rivero.
\newblock The {L}amperti representation of real-valued self-similar {M}arkov
  processes.
\newblock \emph{Bernoulli}, 19\penalty0 (5B):\penalty0 2494--2523, 2013.
\newblock ISSN 1350-7265. \MR{3160562}

\bibitem[Chen et~al.(2007)Chen, Fukushima, and Ying]{CFY-ext}
Z.-Q. Chen, M.~Fukushima, and J.~Ying.
\newblock \emph{Extending {M}arkov processes in weak duality by {P}oisson point
  processes of excursions}, volume~2 of \emph{Abel Symp.}, pages 153--196.
\newblock Springer, Berlin, 2007. \MR{2397787 (2010b:60207)}

\bibitem[Chybiryakov(2006)]{Chy-Lam}
O.~Chybiryakov.
\newblock The {L}amperti correspondence extended to {L}\'evy processes and
  semi-stable {M}arkov processes in locally compact groups.
\newblock \emph{Stochastic Process. Appl.}, 116\penalty0 (5):\penalty0
  857--872, 2006.
\newblock ISSN 0304-4149. \MR{2218339 (2006m:60058)}

\bibitem[Cordero(2010)]{Cor-thesis}
F.~Cordero.
\newblock \emph{On the excursion theory for the symmetric stable L\'evy
  processes with index $\alpha \in ]1,2]$ and some applications}.
\newblock PhD thesis, Universit\'e Pierre et Marie Curie -- Paris VI, 2010.

\bibitem[Fitzsimmons and Getoor(1995)]{FG-bridges}
P.~J. Fitzsimmons and R.~K. Getoor.
\newblock Occupation time distributions for {L}\'evy bridges and excursions.
\newblock \emph{Stochastic Process. Appl.}, 58\penalty0 (1):\penalty0 73--89,
  1995.
\newblock ISSN 0304-4149. \MR{1341555 (96d:60117)}

\bibitem[Gradshteyn and Ryzhik(2007)]{GR}
I.~S. Gradshteyn and I.~M. Ryzhik.
\newblock \emph{Table of integrals, series, and products}.
\newblock Elsevier/Academic Press, Amsterdam, seventh edition, 2007.
\newblock ISBN 978-0-12-373637-6.
\newblock Translated from the Russian. Translation edited and with a preface by
  Alan Jeffrey and Daniel Zwillinger. \MR{2360010 (2008g:00005)}

\bibitem[Hubalek and Kuznetsov(2011)]{ECP1601}
F.~Hubalek and A.~Kuznetsov.
\newblock A convergent series representation for the density of the supremum of
  a stable process.
\newblock \emph{Electron. Commun. Probab.}, 16:\penalty0 84--95, 2011.
\newblock ISSN 1083-589X. \MR{2763530 (2012a:60144)}

\bibitem[Ivanovs(2011)]{Iva-thesis}
J.~Ivanovs.
\newblock \emph{One-sided Markov additive processes and related exit problems}.
\newblock PhD thesis, Universiteit van Amsterdam, 2011.

\bibitem[Kingman(1961)]{Kin-convex}
J.~F.~C. Kingman.
\newblock A convexity property of positive matrices.
\newblock \emph{Quart. J. Math. Oxford Ser. (2)}, 12:\penalty0 283--284, 1961.
\newblock ISSN 0033-5606. \MR{0138632 (25 \#2075)}

\bibitem[Kuznetsov(2011)]{Kuz2011}
A.~Kuznetsov.
\newblock On extrema of stable processes.
\newblock \emph{Ann. Probab.}, 39\penalty0 (3):\penalty0 1027--1060, 2011.
\newblock ISSN 0091-1798. \MR{2789582 (2012e:60139)}

\bibitem[Kuznetsov and Pardo(2013)]{KP-HG}
A.~Kuznetsov and J.~C. Pardo.
\newblock Fluctuations of stable processes and exponential functionals of
  hypergeometric {L}\'evy processes.
\newblock \emph{Acta Appl. Math.}, 123:\penalty0 113--139, 2013.
\newblock ISSN 0167-8019. \MR{3010227}

\bibitem[Kyprianou et~al.(2014)Kyprianou, Pardo, and Watson]{KPW-cens}
A.~E. Kyprianou, J.~C. Pardo, and A.~R. Watson.
\newblock Hitting distributions of {$\alpha$}-stable processes via path
  censoring and self-similarity.
\newblock \emph{Ann. Probab.}, 42\penalty0 (1):\penalty0 398--430, 2014.
\newblock ISSN 0091-1798. \MR{3161489}

\bibitem[Lamperti(1972)]{LampertiT}
J.~Lamperti.
\newblock Semi-stable {M}arkov processes. {I}.
\newblock \emph{Z. Wahrscheinlichkeitstheorie und Verw. Gebiete}, 22:\penalty0
  205--225, 1972. \MR{0307358 (46 \#6478)}

\bibitem[Letemplier and Simon(2013)]{LS-hit}
J.~Letemplier and T.~Simon.
\newblock Unimodality of hitting times for stable processes.
\newblock Preprint, \arxivref{1309.5321v2}{math.PR}, 2013.

\bibitem[Maulik and Zwart(2006)]{MZ-ef}
K.~Maulik and B.~Zwart.
\newblock Tail asymptotics for exponential functionals of {L}\'evy processes.
\newblock \emph{Stochastic Process. Appl.}, 116\penalty0 (2):\penalty0
  156--177, 2006.
\newblock ISSN 0304-4149. \MR{2197972 (2007g:60049)}

\bibitem[Miller(1961)]{Mil-convex}
H.~D. Miller.
\newblock A convexity property in the theory of random variables defined on a
  finite {M}arkov chain.
\newblock \emph{Ann. Math. Statist.}, 32:\penalty0 1260--1270, 1961.
\newblock ISSN 0003-4851. \MR{0126886 (23 \#A4180)}

\bibitem[Pant{\'{\i}}(2013)]{Pan-cond}
H.~Pant{\'{\i}}.
\newblock On {L}\'evy processes conditioned to avoid zero.
\newblock Preprint, \arxivref{1304.3191v1}{math.PR}, 2013.

\bibitem[Peskir(2008)]{Pes-law}
G.~Peskir.
\newblock The law of the hitting times to points by a stable {L}\'evy process
  with no negative jumps.
\newblock \emph{Electron. Commun. Probab.}, 13:\penalty0 653--659, 2008.
\newblock ISSN 1083-589X. \MR{2466193 (2010b:60143)}

\bibitem[Rivero(2007)]{Riv-RE2}
V.~Rivero.
\newblock Recurrent extensions of self-similar {M}arkov processes and
  {C}ram\'er's condition. {II}.
\newblock \emph{Bernoulli}, 13\penalty0 (4):\penalty0 1053--1070, 2007.
\newblock ISSN 1350-7265. \MR{2364226 (2008k:60173)}

\bibitem[Sato(1999)]{Sato}
K.~Sato.
\newblock \emph{L\'evy processes and infinitely divisible distributions},
  volume~68 of \emph{Cambridge Studies in Advanced Mathematics}.
\newblock Cambridge University Press, Cambridge, 1999.
\newblock ISBN 0-521-55302-4. \MR{1739520 (2003b:60064)}

\bibitem[Simon(2011)]{Sim-hit}
T.~Simon.
\newblock Hitting densities for spectrally positive stable processes.
\newblock \emph{Stochastics}, 83\penalty0 (2):\penalty0 203--214, 2011.
\newblock ISSN 1744-2508. \MR{2800088 (2012g:60040)}

\bibitem[Vuolle-Apiala(1994)]{VA-Ito}
J.~Vuolle-Apiala.
\newblock It\^o excursion theory for self-similar {M}arkov processes.
\newblock \emph{Ann. Probab.}, 22\penalty0 (2):\penalty0 546--565, 1994.
\newblock ISSN 0091-1798. \MR{1288123 (95i:60063)}

\bibitem[Yano et~al.(2009{\natexlab{a}})Yano, Yano, and Yor]{YYY-laws}
K.~Yano, Y.~Yano, and M.~Yor.
\newblock On the laws of first hitting times of points for one-dimensional
  symmetric stable {L}\'evy processes.
\newblock In \emph{S\'eminaire de {P}robabilit\'es {XLII}}, volume 1979 of
  \emph{Lecture Notes in Math.}, pages 187--227. Springer, Berlin,
  2009{\natexlab{a}}. \MR{2599211 (2011e:60107)}

\bibitem[Yano et~al.(2009{\natexlab{b}})Yano, Yano, and Yor]{YYY-pen}
K.~Yano, Y.~Yano, and M.~Yor.
\newblock Penalising symmetric stable {L}\'evy paths.
\newblock \emph{J. Math. Soc. Japan}, 61\penalty0 (3):\penalty0 757--798,
  2009{\natexlab{b}}.
\newblock ISSN 0025-5645. \MR{2552915}

\end{thebibliography}

\end{document}